\documentclass[pdflatex,sn-mathphys-num]{sn-jnl}

\usepackage{graphicx}%
\usepackage{multirow}%
\usepackage{amsmath,amssymb,amsfonts}%
\usepackage{amsthm}%
\usepackage{mathrsfs}%
\usepackage{mathtools}
\usepackage[title]{appendix}%
\usepackage{xcolor}%
\usepackage{textcomp}%
\usepackage{manyfoot}%
\usepackage{booktabs}%
\usepackage{algorithm}%
\usepackage{algorithmicx}%
\usepackage{algpseudocode}%
\usepackage{algcompatible}
\usepackage{listings}%
\usepackage{xspace}

\newcommand{\asc}{ASC\xspace}
\newcommand{\fdyn}{\texttt{FedDyn}\xspace}
\newcommand{\scaf}{\texttt{SCAFFOLD}\xspace}
\newcommand{\fdynExt}{\texttt{FedDyn (exact)}\xspace}
\newcommand{\fdynInext}{\texttt{FedDyn (inexact)}\xspace}
\newcommand{\scafSto}{\texttt{SCAFFOLD (stochastic)}\xspace}
\newcommand{\scafFull}{\texttt{SCAFFOLD (full gradient)}\xspace}

\theoremstyle{thmstyleone}%
\newtheorem{theorem}{Theorem}
\newtheorem{lemma}[theorem]{Lemma}
\newtheorem{assumption}[theorem]{Assumption}

\theoremstyle{thmstyletwo}%
\newtheorem{remark}{Remark}%

\theoremstyle{thmstylethree}%
\providecommand{\dataname}{Data Availability}
\def\datablock{\par\addvspace{17pt}\small\rmfamily
\trivlist\if!\dataname!\item[]\else
\item[\hskip\labelsep
{\bfseries\dataname}]\fi}

\newenvironment{data_availability}{\begin{datablock}}{\end{datablock}}

\providecommand{\confname}{Conflicts of Interest}
\def\confblock{\par\addvspace{17pt}\small\rmfamily
\trivlist\if!\confname!\item[]\else
\item[\hskip\labelsep
{\bfseries\confname}]\fi}

\newenvironment{conflicts_of_interest}{\begin{confblock}}{\end{confblock}}

\begin{document}

\title[Federated Optimization under Averaged Strong Convexity Assumption]{Federated Optimization under Averaged Strong Convexity Assumption}


\author[1]{\fnm{Sirong} \sur{Dai}}\email{sirongdai@link.cuhk.edu.cn}

\author*[1]{\fnm{Ming} \sur{Yan}}\email{yanming@cuhk.edu.cn}

\affil[1]{\orgdiv{School of Data Science}, \orgname{The Chinese University of Hong Kong, Shenzhen}, \orgaddress{\city{Shenzhen}, \postcode{518172}, \state{Guangdong}, \country{China}}}


\abstract{This work studies federated optimization under the averaged strong convexity assumption, which allows the global objective to be strongly convex even when the individual local objectives are only convex. Under this assumption, we establish new convergence guarantees for the federated optimization algorithms \fdyn and \scaf under partial client participation. For \fdyn, we derive the communication complexity for both exact and inexact subproblem solutions and further establish the expected parallel gradient complexity for the inexact case. For \scaf, we strengthen the original convergence analysis by establishing a last-iterate error bound that holds from the first communication round, and we derive the corresponding communication complexity. Our analysis establishes almost sure convergence of the iterates and linear convergence of the last iterate in mean square for \fdyn with exact subproblem solutions and for \scaf with exact local gradients. Finally, we illustrate the theoretical results through numerical experiments on distributed least squares and binary logistic regression. }




\maketitle

\section{Introduction}
\label{sec:intro}

Federated learning (FL) is a distributed learning paradigm in which multiple devices jointly train a shared model while keeping raw data on-device~\cite{federatedlearningfedavg}. FL has been deployed in a range of applications, including personalized keyboard prediction on mobile phones~\cite{hard2019federatedlearningmobilekeyboard}, collaborative medical image analysis across hospitals~\cite{adnan2022federated}, recommendation systems spanning multiple e-commerce platforms~\cite{li2024federated}, and learning over IoT/edge networks~\cite{vahabi2025federated}.

FL is typically formulated as a federated optimization problem~\cite{federatedoptimization}, where the server seeks to minimize the average of client-specific objectives,
\[\bar{f}(x) = \frac{1}{N} \sum_{i=1}^{N} f_i(x).\]
Each client $i$ can evaluate and (approximately) optimize its local loss $f_i$, while a central server coordinates the learning process.

Practical FL systems face several well-known challenges, including limited communication, data heterogeneity (non-i.i.d.\ data), privacy constraints, and \emph{partial participation}. Partial participation means that in each communication round only a subset of clients is available to perform local updates and report to the server, due to device availability, intermittent connectivity, or communication-budget limits.

To cope with limited communication and heterogeneity, modern FL methods typically perform multiple local updates per communication round and introduce correction/regularization mechanisms to control client drift. A canonical baseline is \texttt{FedAvg}~\cite{federatedlearningfedavg} (also known as \texttt{local SGD}~\cite{stich2018local}), which reduces communication by letting each client take several local gradient steps before aggregation. Although effective in near i.i.d.\ regimes, \texttt{FedAvg} can converge to a biased solution under heterogeneous data.
To address this, \scaf~\cite{karimireddy2020scaffold} uses control variates to correct client drift and can achieve exact convergence under partial participation with unbiased gradients. From a complementary viewpoint, \fdyn~\cite{acarfederated} is derived from a primal--dual (augmented Lagrangian) perspective and employs a dynamic regularization mechanism that promotes global consistency across clients.

Standard convergence analyses for federated optimization often assume that every local objective $f_i$ is strongly convex. This can be restrictive in practice: some clients may have limited, noisy, or low-rank data, making their local problems flat or underdetermined and hence not strongly convex.

To better reflect such settings, we adopt the \emph{averaged strong convexity} (\asc) assumption. Under \asc, each $f_i$ is only convex, but the global objective $\bar f$ is $\mu$-strongly convex. 

For instance, consider distributed linear regression where client $i$
minimizes $f_i(x) = \frac{1}{2} \|A_i x - b_i\|^2$ with $A_i \in \mathbb{R}^{n_i \times d}$. If $n_i<d$, then $A_i^\top A_i$ is rank-deficient and $f_i$ is not strongly convex. However, if the aggregated data provide full rank, i.e., $\sum_{i=1}^N A_i^\top A_i \succ 0$, then the global objective $\bar f$ is strongly convex. This illustrates the motivation for \asc: global curvature—and hence a unique optimum—can arise from the collective dataset even when each local problem is individually ill-conditioned.

We study federated optimization under the \asc assumption. Our main contributions are as follows:

\begin{enumerate}
    \item \textbf{\fdyn\ under \asc.}
    We provide a new convergence analysis of \fdyn under \asc, allowing for inexact local subproblem solutions. In the exact-solve case, we prove linear convergence in iterate error and almost sure convergence to the unique minimizer, strengthening the guarantees in the original \fdyn analysis.

    \item \textbf{\scaf\ under \asc.}
    We analyze \scaf under the same \asc setting. With exact (noise-free) local gradients, we establish linear convergence in iterate error and almost sure convergence to the unique minimizer, with a rate comparable to \fdyn. Compared to the original analysis of \scaf, our result allows for larger stepsizes and does not require a minimum number of communication rounds.
\end{enumerate}

The remainder is organized as follows. 
We review related work on federated optimization and \asc in Section~\ref{sec:relatedwork}.
Section~\ref{sec:preliminaries} introduces notation, problem formulation, and assumptions.
We develop the main results, proving linear convergence of \fdyn and \scaf under \asc in Section~\ref{sec:theory}.
Section~\ref{sec:experiments} reports experiments on distributed least squares and logistic regression.
Section~\ref{sec:conclusion} concludes and discusses future directions.


\section{Related Work}
\label{sec:relatedwork}

Federated optimization was popularized by the seminal \texttt{FedAvg} method~\cite{federatedlearningfedavg,zhao2018federated,stich2018local,khaled2020tighter,liconvergence}, also referred to as \texttt{local SGD}~\cite{stich2018local,khaled2020tighter,khaled2019first,woodworth2020local}.  Its central idea is to reduce communication by performing multiple local updates on each client between communication rounds. However, under heterogeneous (non-i.i.d.) data, \texttt{FedAvg} typically does not converge to the exact minimizer and instead approaches a bias-dominated neighborhood.
To mitigate this effect, \texttt{FedProx}~\cite{li2020federated} adds a quadratic proximal term to discourage local iterates from drifting too far from the global model. While this regularization can improve stability, it does not fully eliminate heterogeneity-induced bias, and \texttt{FedProx} may still fail to reach the global optimum in highly non-i.i.d.\ regimes.

To better handle heterogeneity, several methods draw on ideas from distributed optimization and operator-splitting. \texttt{FedSplit}~\cite{pathak2020fedsplit} applies Peaceman--Rachford splitting~\cite{peaceman1955numerical} to a reformulated FL objective. \texttt{FedPD}~\cite{zhang2021fedpd} is based on a distributed ADMM formulation~\cite{boyd2011distributed,yan2017self} and provides convergence guarantees with inexact proximal updates. Related \texttt{FedADMM}-type methods have also been proposed to reduce heterogeneity-induced bias. For example,~\cite{song2025fedadmm} studies inexact subproblem solutions with adaptive penalty updates, but focuses on full participation.
To support partial participation, \fdyn~\cite{acarfederated} introduces an ADMM-related method with uniform client sampling. Subsequent work analyzes \texttt{FedADMM} variants under alternative participation models:~\cite{gong2022fedadmm} studies nonconvex objectives with nonuniform (independent) sampling probabilities;~\cite{zhou2023federated} considers a different participation mechanism; and~\cite{wang2023beyond} provides a unified treatment of several \texttt{FedADMM}-type algorithms with nonconvex complexity bounds. For composite objectives with an additional nonsmooth convex regularizer,~\cite{tran2021feddr} applies Douglas-Rachford splitting~\cite{lions1979splitting} to nonconvex federated problems and derives corresponding complexity guarantees.

From a primal (gradient-based) viewpoint, \scaf~\cite{karimireddy2020scaffold} introduces a control-variate mechanism inspired by SAGA~\cite{defazio2014saga}. It treats heterogeneity-induced client drift as a variance effect and corrects it using stored gradient estimates, while still allowing multiple local steps between communication rounds. More recently,~\cite{luo2025revisiting} refines the analysis of \scaf in the nonconvex regime.
Another related line is \texttt{ProxSkip}~\cite{mishchenko2022proxskip}, which uses probabilistic server aggregation and achieves optimal communication complexity via a different analysis. Several works further improve communication efficiency:~\cite{huangxinmengSCALLION} develops unbiased and biased compressed variants of \scaf;~\cite{momentumScaf} combines \scaf with local momentum and variance reduction to improve nonconvex complexity; and~\cite{you2025fedsumfamilyefficientfederated} proposes an algorithm with lower communication cost than \scaf and analyzes it under various partial-participation models.

The \asc assumption relaxes the standard requirement that each local objective $f_i$ be strongly convex. Instead, it only requires the global objective $\bar f$ to be strongly convex, while allowing each $f_i$ to be only convex. This model better reflects practical FL deployments, where individual clients may have limited, biased, or highly heterogeneous data.
In federated optimization, Li et al.~\cite{Extrapolation} introduced
\textsc{FedExProx}, which adds a server-side extrapolation step to
\textsc{FedProx}, and proved linear convergence under \asc and the
additional interpolation assumption that all local objectives share
a common minimizer. Closely related global (averaged) strong convexity assumptions have appeared in decentralized optimization~\cite{shi2015extra,li2021linear}: e.g., \textsc{EXTRA}~\cite{shi2015extra} establishes linear convergence under global strong convexity, and~\cite{li2021linear} provides refined analyses for both \textsc{EXTRA} and \textsc{NIDS} under the same condition. Similar ideas also arise in stochastic optimization; for instance,~\cite{mishchenko2020random} analyzes random reshuffling for finite-sum problems under an assumption closely related to \asc.


\section{Preliminaries}
\label{sec:preliminaries}

This section introduces the notation, formal problem formulation, and assumptions used throughout.

{\bf Notation.} 
Vectors are denoted by lowercase letters (e.g., $x, y \in \mathbb{R}^d$). The Euclidean inner product and norm are $\langle x, y \rangle$ and $\| x \| \coloneqq \sqrt{\langle x, x \rangle}$, respectively.
Given a probability space with a $\sigma$-field $\mathcal{F}$, we write $\mathbb{E}[\,\cdot \mid \mathcal{F}\,]$ for conditional expectation.
We use $[N] \coloneqq \{1, 2, \dots, N\}$ to index the $N$ clients.

We study the standard federated optimization problem
\begin{equation}
    \label{eq:problem}
    \min_{x \in \mathbb{R}^d} \bar{f}(x) \coloneqq \frac{1}{N} \sum_{i=1}^{N} f_i(x),
\end{equation}
where $N$ is the number of clients and each $f_i:\mathbb{R}^d\to\mathbb{R}$ represents the empirical (or population) risk associated with client $i$'s local data.
We consider two common modeling settings:
\begin{itemize}
    \item \textbf{Finite-sum:} 
        Client $i$ holds a dataset $\mathcal D_i$ with $n_i\coloneqq|\mathcal D_i|$ samples, and
    \begin{equation*}
        f_i(x) = \frac{1}{n_i} \sum_{\xi \in \mathcal{D}_i} \ell_i(x; \xi),
    \end{equation*}
    where $\ell_i(\cdot; \xi): \mathbb{R}^d \to \mathbb{R}$ is the per-sample loss.

    \item \textbf{Stochastic (expectation):}     Client $i$ has access to a distribution $\mathcal P_i$, and
    \begin{equation*}
        f_i(x) = \mathbb{E}_{\xi \sim \mathcal{P}_i} \big[ \ell_i(x; \xi) \big].
    \end{equation*}
\end{itemize}

The objective $\bar f$ in \eqref{eq:problem} measures average performance across clients. Our goal is to compute a minimizer $x^*$ of $\bar f$ under
the operational constraints of FL, including limited communication, partial participation, and heterogeneous (non-i.i.d.) local data.

Our analysis relies on two standard structural assumptions.

\begin{assumption}[$L$-Smoothness of Local Objectives]
    \label{ass:smooth}
    Each $f_i:\mathbb{R}^d\to\mathbb{R}$ is continuously differentiable and has an $L$-Lipschitz gradient for some $L>0$ common to all clients; i.e., for all $i\in[N]$ and all $x,y\in\mathbb{R}^d$,
    \begin{equation*}
        \| \nabla f_i({x}) - \nabla f_i({y}) \| \leq L \| {x} - {y} \|.
    \end{equation*}
\end{assumption}

\begin{assumption}[Averaged Strong Convexity (\asc)]
    \label{ass:ASC}
    The global objective $\bar f$ is $\mu$-strongly convex for some $\mu>0$, while each local objective $f_i$ is only convex.
\end{assumption}

Recall that a differentiable function $g:\mathbb{R}^d\to\mathbb{R}$ is $\mu$-strongly convex if, for all $x,y\in\mathbb{R}^d$,
\begin{equation*}
    g({y}) \geq g({x}) + \langle \nabla g({x}), {y} - {x} \rangle + \frac{\mu}{2} \| {y} - {x} \|^2.
\end{equation*}
In particular, $\mu$-strong convexity implies that $g$ admits a unique minimizer.

\begin{remark}
Under Assumptions~\ref{ass:smooth}--\ref{ass:ASC}, $\bar f$ is $L$-smooth and $\mu$-strongly convex, hence has a unique minimizer
$x^*\coloneqq \arg\min_x \bar f(x)$. Moreover, heterogeneity is allowed: in general $\nabla f_i(x^*)\neq 0$ for individual clients.
\end{remark}


\section{Federated Optimization under Averaged Strong Convexity}
\label{sec:theory}

In this section, we present \fdyn and \scaf and analyze the convergence
of each algorithm. Both analyses rely on the following inequality for
smooth convex functions.

\begin{lemma}[Key Inequality for Smooth Convex Functions]
    \label{lem:mainineq}
    Let $f:\mathbb{R}^d\to\mathbb{R}$ be convex and $L$-smooth. Then, for any $x,y,z\in\mathbb{R}^d$,
    \begin{equation}
        \label{eq:mainineq}
        \langle \nabla f(z), x - y \rangle 
        \geq f(x) - f(y) 
        + \frac{1}{2L} {\| \nabla f(z) - \nabla f(y) \|^2} 
        - \frac{L}{2} \| x - z \|^2.
    \end{equation}
\end{lemma}

\noindent\textbf{Proof strategy under \asc.}
When each $f_i$ is strongly convex, we can use its curvature to bound
each client's local iterates before aggregation. Under \asc, strong
convexity applies only after aggregation, so the challenge is to obtain
sufficiently sharp local bounds using only convexity and smoothness.
We use~\eqref{eq:mainineq} to retain a nonnegative gradient-difference
term in each local bound. After aggregation, this term and the strong
convexity of the global objective yield the Lyapunov recursions for
both \fdyn and \scaf. We prove Lemma~\ref{lem:mainineq} in
Appendix~\ref{sec:preliminary-lemmas}.

\subsection{FedDyn}

We begin with \fdyn. To formalize the randomness arising from client sampling, we introduce the $\sigma$-field $\mathcal{F}_r$, which captures all randomness generated up to the beginning of round $r$. Specifically,
\begin{equation*}
    \mathcal{F}_r = \sigma\{S^{0}, S^{1}, \ldots, S^{r-1}\},
\end{equation*}
that is, the $\sigma$-field generated by the sequence of sampled client subsets $\{S^t\}_{t=0}^{r-1}$.
\begin{algorithm}[!h]
\caption{\fdyn}
\label{alg:FedDyn}
\begin{algorithmic}[1]
\small
\REQUIRE \textbf{Server:} initial $x^0=y^0$, $h^0=0$, participation clients $s$, communication rounds $R$.
\REQUIRE \textbf{Clients:} parameter $\gamma>0$, dual variables $\lambda_i^0=0$ for $i\in[N]$.
\FOR{$r=0$ \textbf{to} $R-1$}
    \STATE Sample a subset $S^r\subseteq[N]$ uniformly at random with $|S^r|=s$.
    \STATE Server broadcasts $y^r$ to all clients in $S^r$.
    \STATE \textbf{Client update (in parallel) for each $i\in S^r$:}
        \STATE\hspace{5mm} Approximately solve
        \[
        x_i^{r+1}\approx \arg\min_{x}\left\{
        f_i(x)-\langle \lambda_i^r,\,x-y^r\rangle+\frac{\gamma}{2}\|x-y^r\|^2
        \right\}.
        \]
        \STATE\hspace{5mm} Update $\lambda_i^{r+1}=\lambda_i^r-\gamma(x_i^{r+1}-y^r)$.
        \STATE\hspace{5mm} Send $\bigl(x_i^{r+1},\lambda_i^{r+1}\bigr)$ to the server.
    \STATE \textbf{Inactive clients:} for each $i\notin S^r$, set $x_i^{r+1}=x_i^r$ and $\lambda_i^{r+1}=\lambda_i^r$.
    \STATE \textbf{Server update:}
        \STATE\hspace{5mm} $h^{r+1}=h^r-\frac{\gamma}{N}\sum_{i\in S^r}(x_i^{r+1}-y^r)
        \;\;{\color{blue}\Bigl(=\frac{1}{N}\sum_{i=1}^N\lambda_i^{r+1}\Bigr)}$.
        \STATE\hspace{5mm} $y^{r+1}=\frac{1}{s}\sum_{i\in S^r}x_i^{r+1}-\frac{1}{\gamma}h^{r+1}$.
\ENDFOR
\end{algorithmic}
\end{algorithm}
\fdyn provides a general framework that accommodates partial participation.  When all clients participate in every round (i.e., full participation), \fdyn reduces to \texttt{FedADMM}~\cite{zhang2021fedpd,gong2022fedadmm}.  
When only a single client participates in each iteration, \fdyn becomes \texttt{point-SAGA}~\cite{defazio2016simple,condat2024simple,condat2025stochastic}.

For each client $i \in S^r$, the algorithm solves a subproblem involving its local objective $f_i$.  
Solving these subproblems exactly can be costly, so we allow inexact solutions. Let
\begin{equation*}
    z_i^{r+1} := \arg\min_{x} \left\{
    f_i(x) - \langle \lambda_i^r, x - y^r \rangle + \frac{\gamma}{2} \| x - y^r \|^2
    \right\}
\end{equation*}
denote the exact minimizer, and let $x_i^{r+1}$ be an approximate solution satisfying
\begin{equation*}
    \|x_i^{r+1} - z_i^{r+1}\| \le \epsilon^{r+1},
    \qquad
    \forall\, r+1 \in [R].
\end{equation*}
We use a common tolerance $\epsilon^{r+1}$ for all clients at iteration $r+1$.

We now present the main convergence result of \fdyn under \asc.

\begin{theorem}[{\normalfont\fdyn} under \asc] \label{thm:fdyn}
    Suppose Assumptions~\ref{ass:smooth}-\ref{ass:ASC} hold, and let $x^*$ be the global minimizer of $\bar f$.
    Define the Lyapunov function
    \begin{equation*}
    \psi^r \coloneqq \|x^r-x^*\|^2
    + A_1 \sum_{i=1}^N\|\lambda_i^r-\nabla f_i(x^*)\|^2 .
    \end{equation*}   
    \noindent
    \textbf{(i) Inexact subproblem solution.}
    Fix arbitrary constants $\theta_1>0$ and $\theta_2>0$, and define
\begin{equation*}
A_1 \;\coloneqq\; \frac{1}{s}\Bigl[\frac{1}{(1+\theta_2)\gamma L}-\frac{3+\theta_1}{\gamma^2}\Bigr],
\qquad
A_2 \;\coloneqq\; \frac12\Bigl(1+\frac{1}{\theta_1}\Bigr)
\;+\;\frac{\gamma}{\theta_2 L}
\;+\;\frac{2\gamma}{\mu}.
\end{equation*}
Assume each client solves its subproblem inexactly as specified earlier, and choose a stepsize $\gamma>0$ such that
\begin{equation*}
\frac{1}{\gamma}\;\le\;\frac{2s}{4(3+\theta_1)(1+\theta_2)sL + N\mu}\le\;\frac{1}{6L}.
\end{equation*}
Then the Lyapunov sequence satisfies, for all $r$,
\begin{equation*}
\mathbb{E}\!\left[\psi^{r+1}\mid\mathcal{F}_r\right]
\;\le\;
\Bigl(1-\frac{\mu}{2\gamma}\Bigr)\psi^r
\;+\;
A_2(\epsilon^{r+1})^2.
\end{equation*}
    \smallskip
    \noindent
    \textbf{(ii) Exact subproblem solution.}
    Assume each client solves its subproblem exactly. Define
\begin{equation*}
A_1 \;\coloneqq\; \frac{1}{s}\Bigl(\frac{1}{\gamma L}-\frac{3}{\gamma^2}\Bigr).
\end{equation*}
If the stepsize satisfies
\begin{equation*}
\frac{1}{\gamma} \;\le\; \frac{s}{6sL+N\mu}\le\frac{1}{6L},
\end{equation*}
then the Lyapunov sequence contracts linearly:
\begin{equation*}
\mathbb{E}\!\left[\psi^{r+1}\mid\mathcal{F}_r\right]
\;\le\;
\Bigl(1-\frac{\mu}{\gamma}\Bigr)\psi^r.
\end{equation*}
\end{theorem}
The proof is postponed to Appendix~\ref{app:feddyn}.

Although the exact-subproblem result can be obtained as a special case of the inexact-subproblem analysis (e.g., by setting the inexactness to zero), we state it separately because it yields a strictly faster contraction factor: $1-{\mu\over\gamma}$, compared with $1-{\mu\over2\gamma}$ in the inexact case.

\noindent\textbf{Consequences of the inexact recursion.}
Under Theorem~\ref{thm:fdyn}(i), choose
\(\theta_1=\theta_2=1\) and \(\gamma=16L+N\mu/s\), and let \(A_2=1+\gamma/L+2\gamma/\mu\).
Each participating client solves its local subproblem by Nesterov's
accelerated gradient method~\cite{nesterov2004introductory},
initialized at \(y^{r-1}\).
To achieve \(\mathbb E\|x^R-x^\star\|^2\le\epsilon\), the communication
rounds \(R\) and expected parallel gradient evaluations
\(\mathbb E[G_{\rm par}(R)]\) satisfy the following bounds.

\begin{itemize}
    \item \textbf{Uniform inexactness.}
    For \(0<\epsilon<\psi^0\), set
    \((\epsilon^r)^2=\mu\epsilon/(4\gamma A_2)\) in every round. Then
    \begin{align*}
    R
    &=\mathcal O\!\left(
    \left(\frac L\mu+\frac Ns\right)\log\frac{1}{\epsilon}
    \right),\\
    \mathbb E[G_{\rm par}(R)]
    &=\mathcal O\!\left(
    \left(\frac L\mu+\frac Ns\right)\log\frac{1}{\epsilon}\,
    \log\!\left(2+\frac{sL^2}{\mu^2\epsilon}\right)
    \right).
    \end{align*}

    \item \textbf{Geometrically decaying inexactness.}
    For \(0<\epsilon<1+\psi^0\), set
    \[
    (\epsilon^r)^2
    =\frac{\mu}{8\gamma A_2}
     \left(1-\frac{5\mu}{8\gamma}\right)^{r-1}.
    \]
    Then
    \begin{align*}
    R
    &=\mathcal O\!\left(
    \left(\frac L\mu+\frac Ns\right)\log\frac{1}{\epsilon}
    \right),\\
    \mathbb E[G_{\rm par}(R)]
    &=\mathcal O\!\left(
    \left(\frac L\mu+\frac Ns\right)\log\frac{1}{\epsilon}
    \left[
    \log\!\left(2+\frac{sL^2}{\mu^2}\right)
    +\log\frac{1}{\epsilon}
    \right]\right).
    \end{align*}
\end{itemize}
Both bounds are proved in Appendix~\ref{app:fdyn-inexact-complexity}.

\noindent\textbf{Consequences of the exact recursion.}
Suppose the conditions of Theorem~\ref{thm:fdyn}(ii) hold (exact subproblem solutions) and define \(\rho \coloneqq 1-\frac{\mu}{\gamma}\in(0,1).\) Then, for all $r$,
\begin{equation}\label{eq:fdyn_exact_recursion}
\mathbb{E}[\psi^{r+1}\mid\mathcal{F}_r]\;\le\;\rho\,\psi^r.
\end{equation}

\begin{itemize}
    \item \textbf{Almost sure convergence.}
    Since $\psi^r\ge 0$, \eqref{eq:fdyn_exact_recursion} implies that $\{\psi^r\}$ is a nonnegative supermartingale. Hence, by the supermartingale convergence theorem~\cite{bertsekas2015convex}, $\psi^r$ converges almost surely to a finite random variable.
    Moreover, taking total expectation in \eqref{eq:fdyn_exact_recursion} gives $\mathbb{E}[\psi^r]\le \rho^r\psi^0\to 0$, which implies $\psi^r\to 0$ almost surely. Consequently, $x^r\to x^*$ almost surely.
    
    \item \textbf{Linear convergence in expectation.}
    Unrolling \eqref{eq:fdyn_exact_recursion} yields, for any $R\ge 1$,
    \begin{equation*}
    \mathbb{E}[\psi^{R}] \;\le\; \rho^{R}\psi^0.
    \end{equation*}
    Since $\psi^{R}\ge \|x^{R}-x^*\|^2$, it suffices to enforce $\mathbb{E}[\psi^{R}]\le \epsilon$ to guarantee
    $\mathbb{E}\|x^{R}-x^*\|^2\le \epsilon$.
    In particular, it suffices to choose
    \[
    R \;=\; \mathcal{O} \Bigl(\frac{\gamma}{\mu}\log\frac{1}{\epsilon}\Bigr).
    \]
    With the smallest admissible choice $\gamma_0\coloneqq \frac{6sL+N\mu}{s}$, this becomes
    \[
    R \;=\;
    \mathcal{O} \Bigl(\bigl(\frac{N}{s}+\frac{L}{\mu}\bigr)\log\frac{1}{\epsilon}\Bigr).
    \]
\end{itemize}

\begin{remark}
Under the stronger assumption that each local objective function $f_i$ is $\mu$-strongly convex and $L$-smooth, the original \fdyn analysis~\cite{acarfederated} established the same iteration complexity, but did not provide a result on a strong (iterate-wise) convergence guarantee. 
More recently,~\cite{condat2024simple} studied a \texttt{Point-SAGA} variant that is equivalent to \fdyn and proved the improved complexity
\begin{equation*}
    \mathcal{O} \Bigl( \bigl( \frac{N}{s} + \sqrt{\frac{N}{s} \cdot \frac{L}{\mu}} \bigr)
    \log \bigl( \frac{1}{\epsilon} \bigr) \Bigr),
\end{equation*}
which exhibits a better dependence on the condition number $\tfrac{L}{\mu}$.  
When $s=1$ (i.e., sampling one client per iteration),  it matches the lower bound in~\cite{woodworth2016tight}.

In contrast, under the weaker assumptions adopted in this work--namely, \asc together with $L$-smoothness--no matching lower bound is currently available.  
To the best of our knowledge, our result is the first to establish linear convergence together with an explicit iteration complexity under this weaker setting for this class of algorithms.  
Whether \fdyn can also attain an improved dependence on the condition number under \asc remains an interesting direction for future research.
\end{remark}


\subsection{SCAFFOLD}

\newcommand{\I}{\hspace{5mm}}
\newcommand{\II}{\hspace*{10mm}}
\begin{algorithm}[!ht]
\caption{\scaf}
\label{alg:Scaffold}
\begin{algorithmic}[1]
\small
\REQUIRE \textbf{Server:} initial model $y^0$, $h^0=0$, global stepsize $\eta_g$, participation clients $s$, communication rounds $R$, local steps $K$.
\REQUIRE \textbf{Clients:} local stepsize $\eta_l$, control variates $\lambda_i^0=0$ for all $i\in[N]$.
\FOR{$r=0$ \textbf{to} $R-1$}
    \STATE Sample a subset $S^r\subseteq[N]$ uniformly at random with $|S^r|=s$.
    \STATE Server broadcasts $(y^r,h^r)$ to all clients in $S^r$.
    \STATE \textbf{Client update (in parallel) for each $i\in S^r$:}
        \STATE\hspace{5mm} Initialize $x_i^{r,0}=y^r$.
        \STATE \I \textbf{Local steps:}
        \FOR{$k=0$ \textbf{to} $K-1$}
            \STATE\hspace{5mm} Sample a stochastic gradient $g_i(x_i^{r,k})$.
            \STATE\hspace{5mm} $x_i^{r,k+1}=x_i^{r,k}-\eta_l\bigl(g_i(x_i^{r,k})-\lambda_i^r+h^r\bigr)$.
        \ENDFOR
        \STATE\hspace{5mm} Update
        \[
        \lambda_i^{r+1}
        =
        \lambda_i^r-h^r+\frac{1}{K\eta_l}\bigl(y^r-x_i^{r,K}\bigr)
        \;\;{\color{blue}
        \Bigl(=\frac{1}{K}\sum_{k=0}^{K-1} g_i(x_i^{r,k})\Bigr)}.
        \]
        \STATE\hspace{5mm} Send $\Delta x_i^r\coloneqq x_i^{r,K}-y^r$ and $\Delta\lambda_i^r\coloneqq \lambda_i^{r+1}-\lambda_i^r$ to the server.
    \STATE \textbf{Inactive clients:} for each $i\notin S^r$, set $\lambda_i^{r+1}=\lambda_i^r$.
    \STATE \textbf{Server update:}
        \STATE\hspace{5mm} $\Delta y^r=\frac{1}{s}\sum_{i\in S^r}\Delta x_i^r$,\quad
        $\Delta h^r=\frac{1}{s}\sum_{i\in S^r}\Delta\lambda_i^r$.
        \STATE\hspace{5mm} $y^{r+1}=y^r+\eta_g\,\Delta y^r$,\quad
        $h^{r+1}=h^r+\frac{s}{N}\Delta h^r
        \;\;{\color{blue}\Bigl(=\frac{1}{N}\sum_{i=1}^N \lambda_i^{r+1}\Bigr)}$.
\ENDFOR
\end{algorithmic}
\end{algorithm}

Since \scaf employs stochastic local updates, we impose the following standard assumption on the stochastic gradients.

\begin{assumption}\label{ass:BoundStoGrad}
At each local update, client $i$ can query a stochastic gradient oracle that returns $g_i(x)$ at any point $x\in\mathbb{R}^d$ such that
\begin{enumerate}
    \item \textbf{Unbiasedness:} $\mathbb{E}\!\left[g_i(x)\mid x\right]=\nabla f_i(x)$;
    \item \textbf{Bounded conditional variance:} $\mathbb{E}\!\left[\|g_i(x)-\nabla f_i(x)\|^2 \mid x\right]\le \sigma^2$.
\end{enumerate}
\end{assumption}

To formalize the randomness due to client sampling and stochastic gradients, we introduce the $\sigma$-algebra generated by all iterates and sampled client sets up to
communication round $r$ as 
\[
\tilde{\mathcal{F}}_r
\;=\;
\sigma\Bigl(
\{x_i^{t,k}\}_{i\in S^{t},\,k=0,\ldots,K},\, S^{t}
:\; t=0,\ldots,r-1
\Bigr).
\]
We now state the convergence guarantee of \scaf under \asc.

\begin{theorem}[{\normalfont\scaf} under \asc]\label{thm:scaf}
Suppose Assumptions~\ref{ass:smooth}-\ref{ass:ASC}, and \ref{ass:BoundStoGrad} hold. Let parameter $\eta_g\ge 1$ and
\[
\eta_l \le \sqrt{\frac{\mu}{L}} \frac{1}{2L(1+K)},
\quad
\tilde{\eta}\coloneqq K\eta_g\eta_l
\le
\tilde{\eta}_{\max}
\coloneqq
\min\Bigl\{\frac{1}{9L}, \frac{2s}{8\frac{1+K}{K}sL+5N\mu}\Bigr\}.
\]
Define the Lyapunov function
\[
\phi^r\coloneqq
\|y^r-y^*\|^2+
\frac {5\tilde{\eta}^2}{2s}\sum_{i=1}^N\|\lambda_i^r-\nabla f_i(y^*)\|^2,
\]
where $y^*$ denotes the global minimizer of $\bar f$. Then, for all $r\ge 0$,
\[
\mathbb{E}\!\left[\phi^{r+1}\mid\tilde{\mathcal{F}}_r\right]
\;\le\;
\Bigl(1-\frac{\tilde{\eta}\mu}{2}\Bigr)\phi^r
\;+\;
\tilde{\eta}\bigl(L\tilde{\eta}^2+9\tilde{\eta}\bigr)\frac{\sigma^2}{K}.
\]

\end{theorem}
The proof is deferred to Appendix~\ref{app:scaffold}.

\noindent\textbf{Consequences of the recursion in Theorem~\ref{thm:scaf}.}
Taking the total expectation of the one-step recursion in Theorem~\ref{thm:scaf} and unrolling it yields, for all $r\ge 1$,
\begin{equation*}
    \begin{aligned}
        \mathbb{E}[\phi^r]
        &\le
        \Bigl(1-\frac{\tilde{\eta}\mu}{2}\Bigr)^r \phi^0 + \tilde{\eta}\bigl(L\tilde{\eta}^2+9\tilde{\eta}\bigr)\frac{\sigma^2}{K}\sum_{t=0}^{r-1}\Bigl(1-\frac{\tilde{\eta}\mu}{2}\Bigr)^t \\
        &\le\;
        \exp\bigl(-\frac{\tilde{\eta}\mu r}{2}\bigr) \phi^0 + \frac 2\mu\bigl(L\tilde{\eta}^2+9\tilde{\eta}\bigr)\frac{\sigma^2}{K}.
    \end{aligned}
\end{equation*}
To mirror the convergence analysis established in~\cite{karimireddy2020scaffold}, we characterize the behavior of the algorithm based on the number of communication rounds $R$. Specifically, when $R$ is small, the algorithm exhibits exponential convergence to a neighborhood, whereas for a large $R$, the convergence becomes sublinear. 
To achieve this, we set the stepsize as $\tilde{\eta}=\min\bigl\{\tilde{\eta}_{\max}, \frac{2}{\mu R}\log(\max\{\mathrm{e},\tfrac{\Phi^0\mu^2R^2}{(\frac{8L}{\mu}+36R)\frac{\sigma^2}{K}}\})\bigr\}$, where $\mathrm{e}$ is the base of the natural logarithm and $\Phi^0\coloneqq\|y^0-y^*\|^2+\frac {5\tilde{\eta}_{\max}^2}{2s}\sum_{i=1}^N\|\nabla f_i(y^*)\|^2$. This naturally leads to two distinct regimes:
\begin{itemize}
    \item \textbf{Case 1 (Small $R$):} When $R$ is small, we have $\tilde{\eta} = \tilde{\eta}_{\max}$. In this regime, the algorithm converges exponentially to a neighborhood of the optimum:
    \begin{align*}
        \mathbb{E}[\phi^R]
        &\leq \tilde{{O}}\biggl(\exp\Bigl(-\min\bigl\{\tfrac{\mu}{32L}, \tfrac{s}{10N}\bigr\}R\Bigr)\Phi^0 + \Bigl(\tfrac{1}{\mu L} + \tfrac{s}{N\mu^2}\Bigr)\tfrac{\sigma^2}{K}\biggr).
    \end{align*}

    \item \textbf{Case 2 (Large $R$):} When $R$ is sufficiently large, the stepsize is determined by  $\tilde{\eta} = \frac{2}{\mu R}\log(\max\{\mathrm{e}, \tfrac{\Phi^0\mu^2R^2}{(\frac{8L}{\mu}+36R)\frac{\sigma^2}{K}}\})$. The convergence becomes sublinear:
    \begin{align*}
        \mathbb{E}[\phi^R]
        &\leq \exp\biggl(-\log(\max\{\mathrm{e},\frac{\Phi^0\mu^2R^2}{(\frac{8L}{\mu}+36R)\frac{\sigma^2}{K}}\})\biggr)\Phi^0 + \tilde{\mathcal{O}}\biggl(\Bigl(\tfrac{L}{\mu^3R^2} + \tfrac{1}{\mu^2R}\Bigr)\tfrac{\sigma^2}{K}\biggr) \\
        &\leq \tilde{\mathcal{O}}\biggl(\Bigl(\tfrac{L}{\mu^3R^2} + \tfrac{1}{\mu^2R}\Bigr)\tfrac{\sigma^2}{K}\biggr),
    \end{align*}
        where $\tilde{\mathcal{O}}(\cdot)$ hides logarithmic factors in the problem parameters.
\end{itemize}
Combining the bounds from both cases, we obtain an overall complexity bound for the iterate error:
\begin{equation*}
    \mathbb{E}\bigl[\|y^R-y^*\|^2\bigr]
    \le
    \tilde{\mathcal{O}}\biggl(
        \exp\Bigl(-\min\bigl\{\tfrac{\mu}{32L}, \tfrac{s}{10N}\bigr\}R\Bigr)\Phi^0
        + \Bigl(\tfrac{L}{\mu^3R^2}+\tfrac{1}{\mu^2R}\Bigr)\tfrac{\sigma^2}{K}
    \biggr).
\end{equation*}

\noindent\textbf{Dependence of the variance term on client participation.}
The recursion in Theorem~\ref{thm:scaf} allows an explicit trade-off between the
variance bound and the contraction of the initial error. For a prescribed
number of communication rounds \(R\geq 1\), choose the constant stepsizes
\begin{equation*}
    \eta_g=1,\qquad
    \eta_\ell=\frac{\widetilde\eta}{K},\qquad
    \widetilde\eta
    =\min\!\left\{
        \widetilde\eta_{\max},\,
        \frac{K}{2L(1+K)}\sqrt{\frac{\mu}{L}},\,
        \frac{2(1+\log R)}{\mu sR}
    \right\}.
\end{equation*}
This choice satisfies the stepsize conditions of Theorem~\ref{thm:scaf}.
Since \(L\widetilde\eta\leq 1/9\), the preceding bound gives
\begin{equation}\label{eq:ScafLinearSpdup}
        \mathbb{E}\|y^R-y^\ast\|^2
        \leq
        \exp\!\left(-\frac{\mu\widetilde\eta R}{2}\right)\phi^0
        +\frac{328}{9}\,
         \frac{\sigma^2(1+\log R)}{\mu^2 sKR}.
\end{equation}
For fixed \(R\) and \(K\), this stepsize choice yields a linear speedup
in the variance bound with respect to \(s\). But this gain comes from 
reducing the effective stepsize. When \(\widetilde\eta=2(1+\log R)/(\mu sR)\), the initial-error
factor is \(\exp(-(1+\log R)/s)\), which decays more slowly as \(s\)
increases. 

\begin{remark}
The original convergence guarantee for \scaf~\cite{karimireddy2020scaffold} is established under the assumption that each local objective function $f_i$ is $\mu$-strongly convex. In this setting, the authors derive an ergodic function-value bound for the averaged iterate $\bar y^R$:

\begin{equation*}
    \small
    \mathbb{E}\,[\bar f(\bar y^R)] - \bar f(y^*)
    \le
    \tilde{\mathcal{O}}\Bigl(
        \tfrac{\sigma^2}{\mu R K s}\Bigl(1+\tfrac{s}{\eta_g^2}\Bigr)
        + \tfrac{N\mu}{s}\,
        \exp\Bigl(
            -\min\Bigl\{\tfrac{s}{30N},\,\tfrac{\mu}{162L}\Bigr\} R
        \Bigr)\hat\phi^0
    \Bigr).
\end{equation*}
    This bound holds provided that the following conditions are met:
\begin{equation*}
    R \;\ge\; \max \Bigl\{\tfrac{162L}{\mu},\,\tfrac{30N}{s}\Bigr\},
    \qquad
    \eta_g \ge 1,
    \qquad
    \tilde{\eta} \le \min\Bigl\{\tfrac{1}{81L},\,\tfrac{s}{15N\mu}\Bigr\}.
\end{equation*}
Note, while iterate convergence is not explicitly derived in~\cite{karimireddy2020scaffold}, a comparable complexity bound for the iterates can be obtained by applying the techniques from~\cite{stich2019unified}. In the original SCAFFOLD bound, setting \(\eta_g=\sqrt{s}\) yields a linear speedup in the variance term with respect to the number 
\(s\) of participating clients, matching the \(1/s\) dependence in~\eqref{eq:ScafLinearSpdup}. Under full participation, the authors of~\cite{scaffoldlinearspeedup} establish 
a linear speedup with respect to the total number of clients \(N\) under the additional assumption that each local objective has a uniformly bounded third derivative.

Our result achieves a leading stochastic term of the same order with respect to $R$ while permitting larger stepsizes. Moreover, our bound holds uniformly for all $R \ge 1$, entirely eliminating the requirement for a minimum number of communication rounds.

\noindent\textbf{Consequences of the exact-gradient recursion.}
If each local step uses the exact gradient, i.e., $g_i(x)=\nabla f_i(x)$, then the recursion reduces to the linear contraction
\begin{equation*}
    \mathbb{E}[\phi^{r+1}\mid\tilde{\mathcal{F}}_r]
    \leq
    \Bigl(1-\frac{\tilde{\eta}\mu}{2}\Bigr)\phi^r.
\end{equation*}

\begin{itemize}
    \item \textbf{Almost sure convergence.}
    By the supermartingale convergence theorem~\cite{bertsekas2015convex}, we have $\phi^r\to 0$ almost surely. Consequently, the primal iterates satisfy $y^r\to y^*$ almost surely.
    
    \item \textbf{Linear convergence in expectation.}
    Taking total expectation and unrolling the recursion gives
    \begin{equation*}
        \mathbb{E}[\phi^R]\leq
        \Bigl(1-\frac{\tilde{\eta}\mu}{2}\Bigr)^R\phi^0.
    \end{equation*}
    Since $\phi^R \ge \|y^R-y^*\|^2$, it suffices to enforce $\mathbb{E}[\phi^R]\le\epsilon$ to guarantee
    $\mathbb{E}[\|y^R-y^*\|^2]\le\epsilon$.
    With $\tilde{\eta}=\tilde{\eta}_{\max}$, the resulting iteration complexity is
    \[
    \mathcal{O}\,\Bigl(\bigl(\frac{N}{s}+\frac{L}{\mu}\bigr)\log\bigl(\frac{1}{\epsilon}\bigr)\Bigr),
    \]
    which matches the complexity of \fdyn with exact subproblem solutions.
\end{itemize}

\end{remark}


\section{Numerical Experiments}
\label{sec:experiments}

This section evaluates the empirical performance of \fdyn and \scaf under the \asc condition. We consider (i) a synthetic distributed least-squares problem and (ii) an $\ell_2$-regularized distributed logistic regression task applied on the MNIST Dataset~\cite{lecun2010mnist}. 
All experiments are implemented in MATLAB R2023a and run on a machine with an Intel(R) Core(TM) i7-8550U CPU.

\subsection{Distributed Least-Squares Problem}

We construct a distributed least-squares problem with $M=1000$ samples
and $d=100$ features. We generate a matrix $A\in\mathbb{R}^{M\times d}$
with i.i.d.\ standard Gaussian entries and set
\[
b=A x_{\mathrm{true}}+0.1\,\xi,
\qquad
x_{\mathrm{true}}\sim\mathcal{N}(0,I_d),
\qquad
\xi\sim\mathcal{N}(0,I_M),
\]
where $A$, $x_{\mathrm{true}}$, and $\xi$ are independent. We compute
the reference solution $x^\star$ from the pooled data.

To vary the number of clients while keeping the global objective fixed,
we randomly permute the $M$ samples once. For each $N\in\{20,60,100\}$,
we divide this ordering into $N$ disjoint blocks and keep the partition
fixed across trials. Each client receives 50 samples for $N=20$ and
10 for $N=100$. For $N=60$, 40 clients receive 17 samples each and the
remaining 20 receive 16 each. Let $(A_i,b_i)$ denote the block of $q_i$
samples assigned to client $i$. We define
\[
f_i(x)=\frac{N}{2M}\|A_i x-b_i\|^2,
\qquad
\bar f(x)=\frac{1}{N}\sum_{i=1}^N f_i(x)
         =\frac{1}{2M}\|Ax-b\|^2.
\]
Thus all three partitions have the same global objective and minimizer.
Since $q_i<d$, every local Hessian is singular, while the global Hessian
$A^\top A/M$ is positive definite almost surely. The observation noise
also ensures that the local objectives almost surely have no common
minimizer.

We compare four variants.
\begin{itemize}
    \item \fdynExt{} solves each local quadratic subproblem exactly.
    \item \fdynInext{} uses Nesterov's accelerated gradient method for
    strongly convex functions, initialized at the current server iterate.
    In round $r+1$, starting at $r=0$, the solver stops when the
    subproblem gradient norm is at most $\gamma\epsilon_{r+1}$, where
    $\epsilon_{r+1}=1.5\times10^{-10}(0.93)^r$.
    \item \scafFull{} uses the full local gradient.
    \item \scafSto{} samples one local observation uniformly at each
    step and scales its gradient by $Nq_i/M$ to obtain an unbiased
    estimate of $\nabla f_i$.
\end{itemize}
Both SCAFFOLD variants use $K=5$ local steps and a global stepsize
$\eta_g=1$.

Although the global Hessian is fixed, the maximum local smoothness
constant $L=\max_i (N/M)\|A_i\|_2^2$ increases from $5.93$ to $14.43$
and $19.31$ for $N=20,60,100$, respectively. We therefore tune the
FedDyn penalty $\gamma$ and the SCAFFOLD local stepsize $\eta$ separately
for each $N$ under full participation ($s=N$). For FedDyn, we select
the candidate that gives the smallest squared iterate error for
\fdynExt{} after 50 rounds. For SCAFFOLD, we start with the best
candidate for \scafFull{} at round 50 and decrease the stepsize along
the grid, if needed, until the mean squared iterate error of \scafSto{}
at round 50 is no greater than its initial value. This yields
$\gamma=2,2,4$ for $N=20,60,100$, respectively, and $\eta=0.01$ for
all three values of $N$. For each $N$, we use the selected parameters
for both variants of each method at all participation rates.

\begin{figure}[!htb]
    \centering

    \includegraphics[trim={95pt 265pt 100pt 270pt},clip,width=0.32\linewidth]{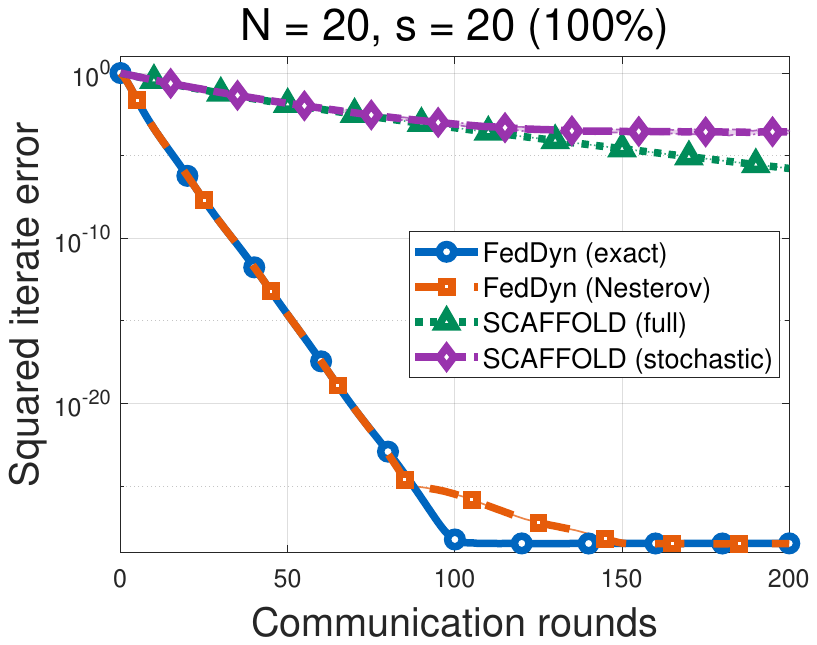}\hfill
    \includegraphics[trim={95pt 265pt 100pt 270pt},clip,width=0.32\linewidth]{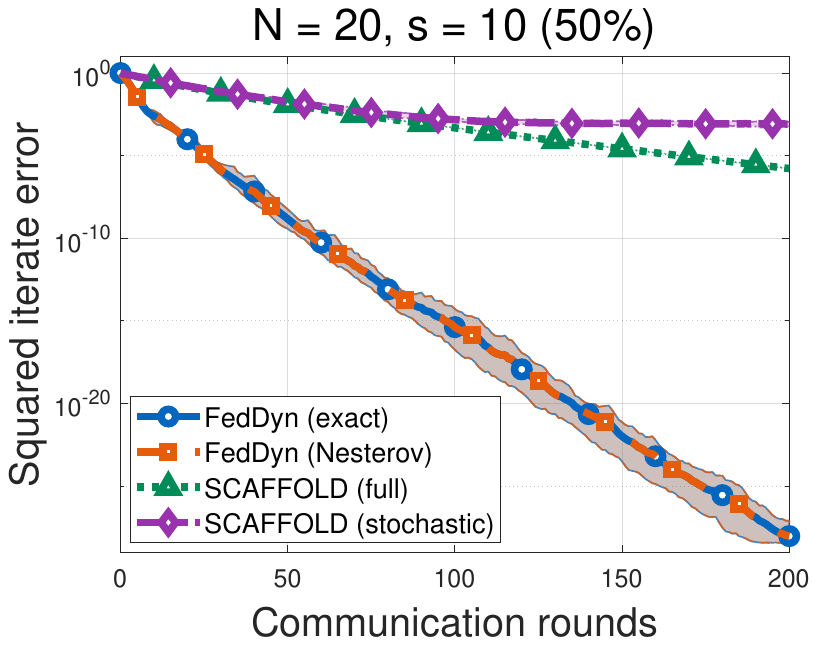}\hfill
    \includegraphics[trim={95pt 265pt 100pt 270pt},clip,width=0.32\linewidth]{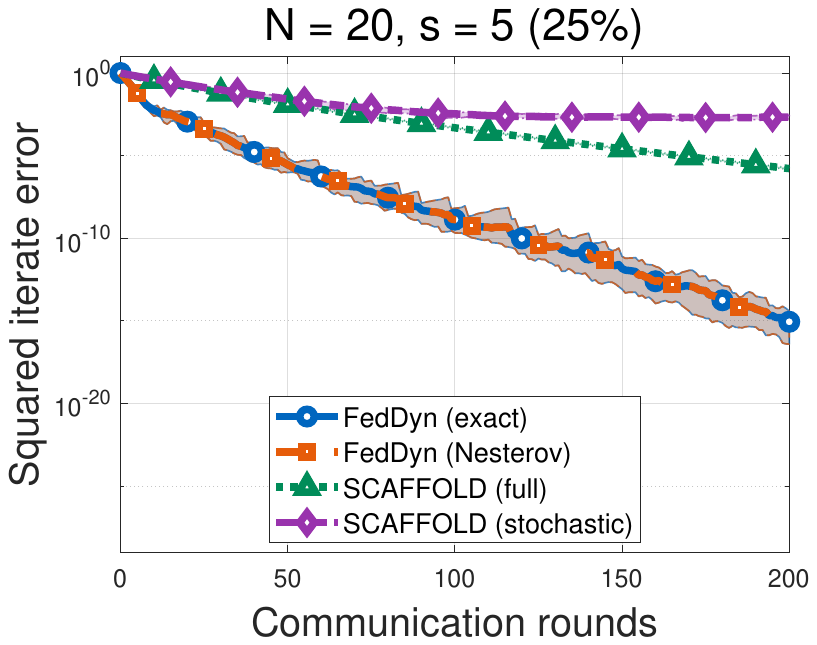}
    \par\medskip

    \includegraphics[trim={95pt 265pt 100pt 270pt},clip,width=0.32\linewidth]{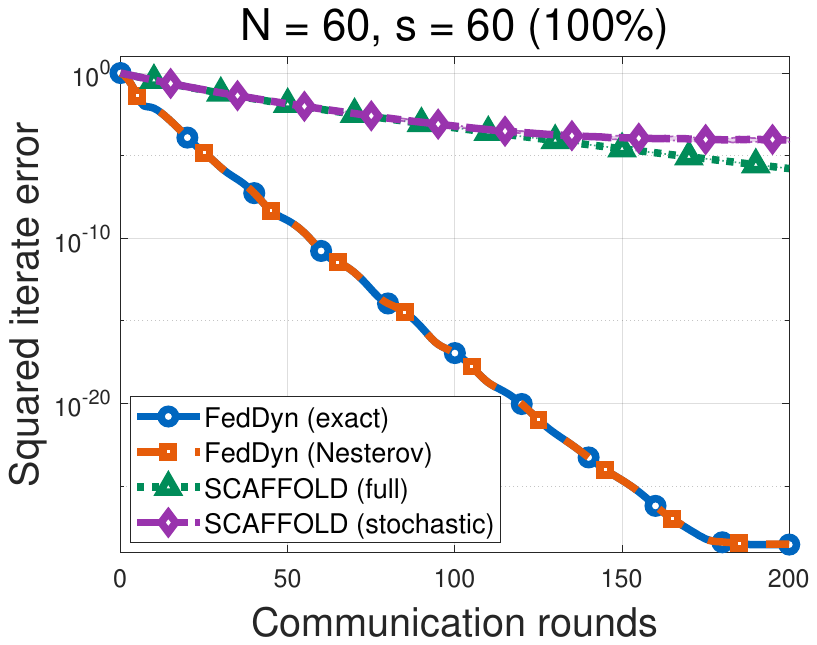}\hfill
    \includegraphics[trim={95pt 265pt 100pt 270pt},clip,width=0.32\linewidth]{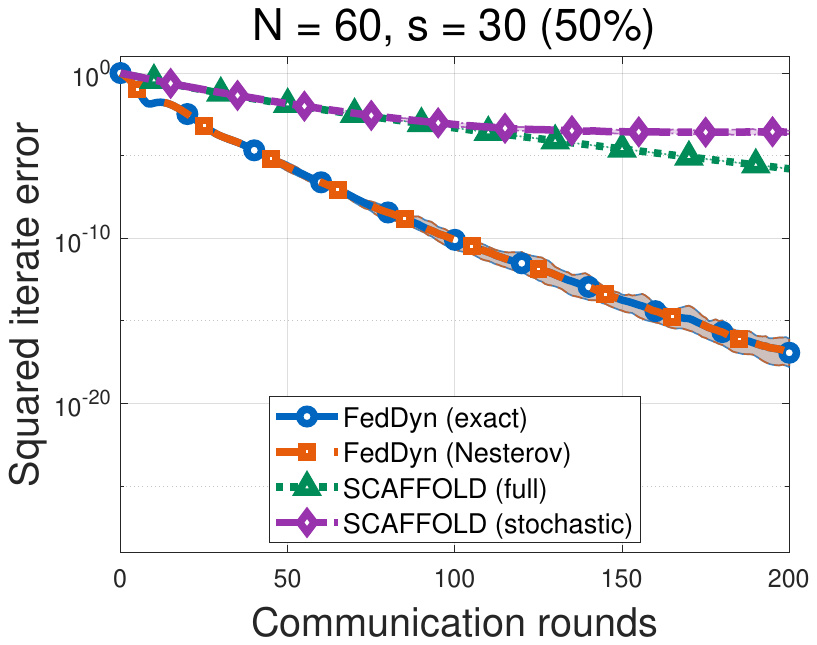}\hfill
    \includegraphics[trim={95pt 265pt 100pt 270pt},clip,width=0.32\linewidth]{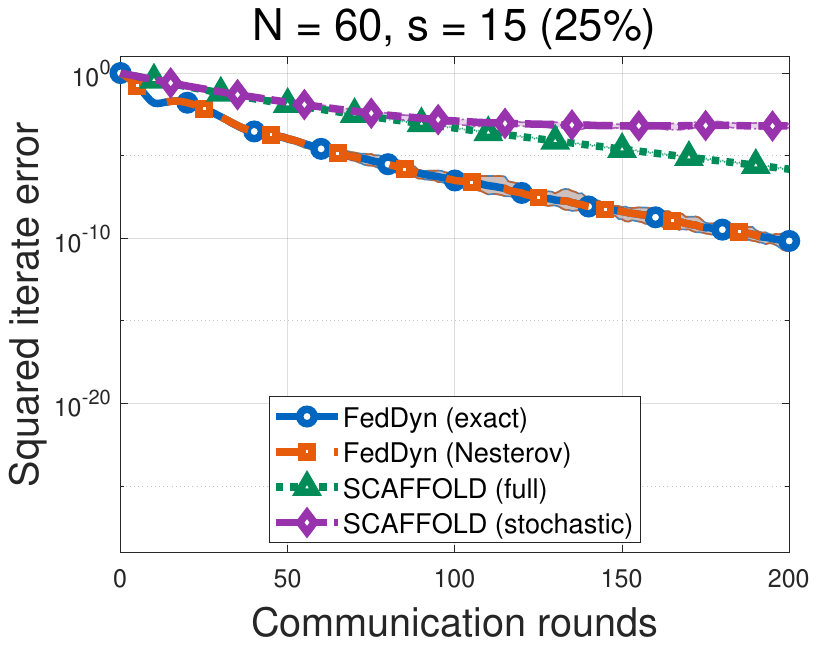}
    \par\medskip

    \includegraphics[trim={95pt 265pt 100pt 270pt},clip,width=0.32\linewidth]{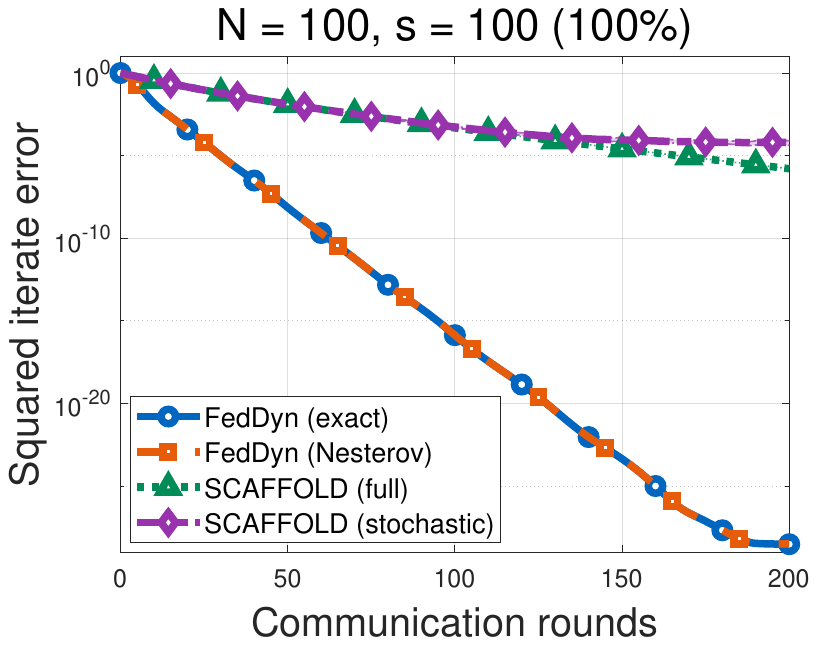}\hfill
    \includegraphics[trim={95pt 265pt 100pt 270pt},clip,width=0.32\linewidth]{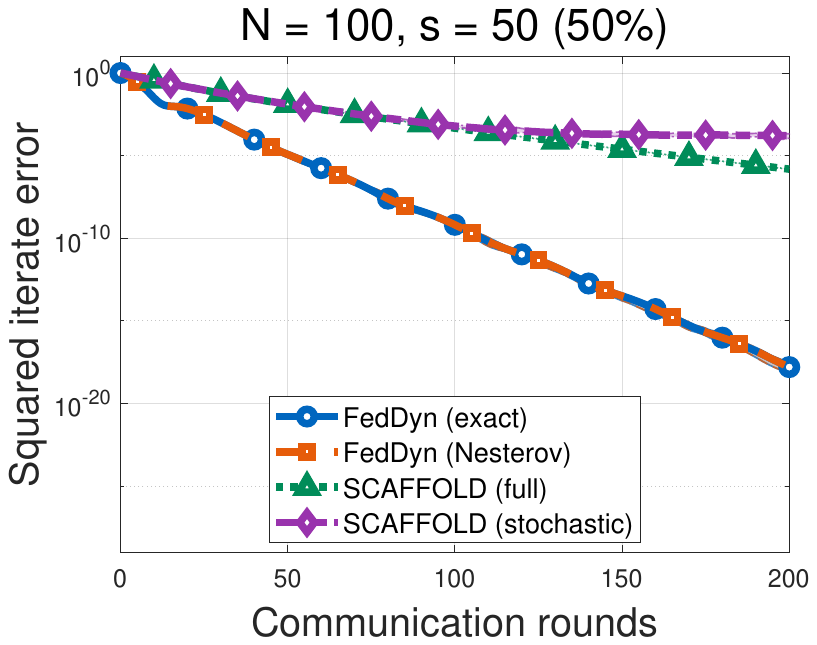}\hfill
    \includegraphics[trim={95pt 265pt 100pt 270pt},clip,width=0.32\linewidth]{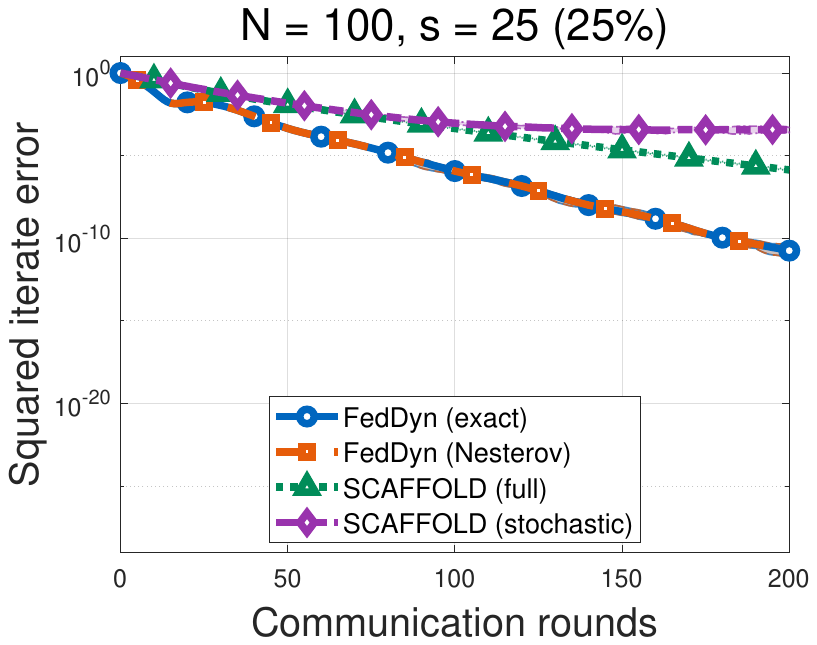}

    \caption{Squared iterate error versus communication rounds for distributed
    least squares under averaged strong convexity. From top to bottom, the
    rows correspond to \(N=20,60,100\) clients. From left to right, the
    columns correspond to \(s=N\), \(s=N/2\), and \(s=N/4\) participating
    clients per round. Curves show the mean squared iterate error over 20 trials, and shaded bands span the minimum and maximum errors observed at each round.}
    \label{fig:convergence_curve_1}
\end{figure}

All methods start from $x^0=x^\star+\mathbf{1}/\sqrt{d}$, so the
initial squared iterate error is one. For each $N$, we consider
$s\in\{N,N/2,N/4\}$ participating clients per round. We sample the
$s$ clients uniformly without replacement and use the same participating
clients for all methods within each trial.
Figure~\ref{fig:convergence_curve_1} shows the mean squared iterate
error over $T=20$ trials for $R=200$ communication rounds. The shaded
bands span the minimum and maximum errors across trials at each round.
To quantify the effect of $N$ on stochastic SCAFFOLD under full
participation, we also report the mean squared iterate error averaged
over the final 50 rounds.

Figure~\ref{fig:convergence_curve_1} shows that, for each fixed $N$,
FedDyn converges faster as the participation rate increases. Under
partial participation, the exact and inexact variants follow similar
trajectories. The squared iterate error of full-gradient SCAFFOLD
continues to decrease over the observed rounds, while that of
stochastic SCAFFOLD levels off. Under full participation, increasing
$N$ from 20 to 60 and 100 reduces the stochastic SCAFFOLD error,
averaged over the final 50 rounds, by factors of 2.76 and 4.00,
respectively. Both reductions fall below the factors of 3 and 5
corresponding to $1/N$ scaling.

\subsection{Distributed Logistic Regression Problem}\label{sec:logisticExperiment}

We consider binary logistic regression on MNIST digits ``4'' and ``9'',
encoded as $0$ and $1$, respectively. We scale pixel values to $[0,1]$
and add a bias coordinate, giving $d=785$. The dataset contains 5,842
images of digit ``4'' and 5,949 images of digit ``9''. For $N=20$, we
retain the original partition with ten clients per digit. For $N=60$
and $N=100$, we shuffle the samples within each class once and divide
the same ordering into $N/2$ disjoint blocks of nearly equal size.
Each client holds one block from a single class, and each partition
is fixed across trials.

To preserve the original global objective, we retain the weight of
each sample under the 20-client partition. Specifically, if sample $j$
belongs to an original block of size $n_b$, we set $w_j=1/(20n_b)$.
Let $D_i$ denote the samples assigned to client $i$ in the partition
with $N$ clients. Define
\[
    \ell_j(x)=\log(1+e^{a_j^\top x})-y_j a_j^\top x,
    \qquad
    f_i(x)=N\sum_{j\in D_i}w_j\ell_j(x)
           +\frac{\rho_i}{2}\|x\|^2.
\]
We set $\rho_i=20\mu$ for $N/20$ clients holding digit ``4'' and
$\rho_i=0$ for all other clients, with $\mu=10^{-2}$. Thus, for every
$N\in\{20,60,100\}$, the global objective is
\[
    \bar f(x)=\frac{1}{N}\sum_{i=1}^N f_i(x)
             =\sum_j w_j\ell_j(x)+\frac{\mu}{2}\|x\|^2.
\]
The three partitions therefore have the same $\mu$-strongly convex
global objective and the same minimizer $x^\star$, while the
unregularized local objectives are convex and lack strong convexity.

We compare \fdynInext{}, \scafFull{}, and \scafSto{}. All methods start
from $x^0=0$ and run for 200 communication rounds. In each configuration,
we fix $s\in\{N,N/2,N/4\}$ and sample $s$ clients uniformly without
replacement at each round. All methods use the same participating
clients within each trial.

As in the least-squares experiment, \fdynInext{} uses Nesterov's
accelerated gradient method for strongly convex subproblems, starting
from $y^r$ and stopping when
$\|\nabla Q_i^r(u)\|\leq\gamma\epsilon_{r+1}$, where
$\epsilon_{r+1}=1.5\times10^{-10}(0.93)^r$. Both SCAFFOLD variants use
$K=5$ local steps and $\eta_g=1$. At each step, \scafFull{} uses the
full local gradient, while \scafSto{} samples $j\in D_i$ uniformly and
uses the unbiased gradient $N\cdot|D_i|\cdot w_j\nabla\ell_j(x)+\rho_i x$.

For each $N$, we tune $\gamma$ for \fdynInext{} and $\eta$ for
\scafFull{} by minimizing the squared iterate error at round 50 under
full participation. We reduce $\eta$ along the candidate grid, if
needed, until the round-50 mean squared error of \scafSto{} is no
greater than its initial value. Both SCAFFOLD variants use this $\eta$,
and the selected parameters are fixed across participation rates.

\begin{figure}[!htb]
    \centering

    \includegraphics[trim={95pt 265pt 100pt 270pt},clip,width=0.32\linewidth]{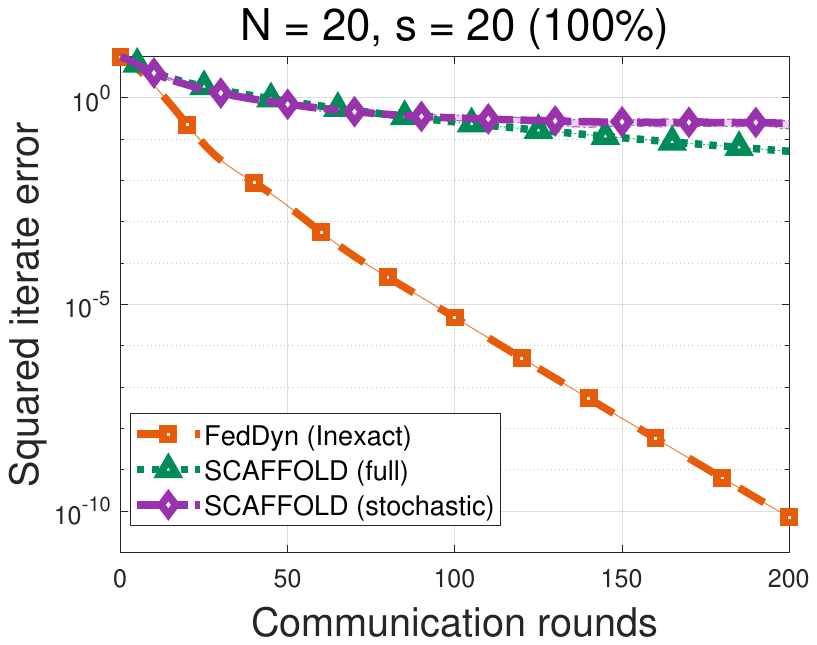}\hfill
    \includegraphics[trim={95pt 265pt 100pt 270pt},clip,width=0.32\linewidth]{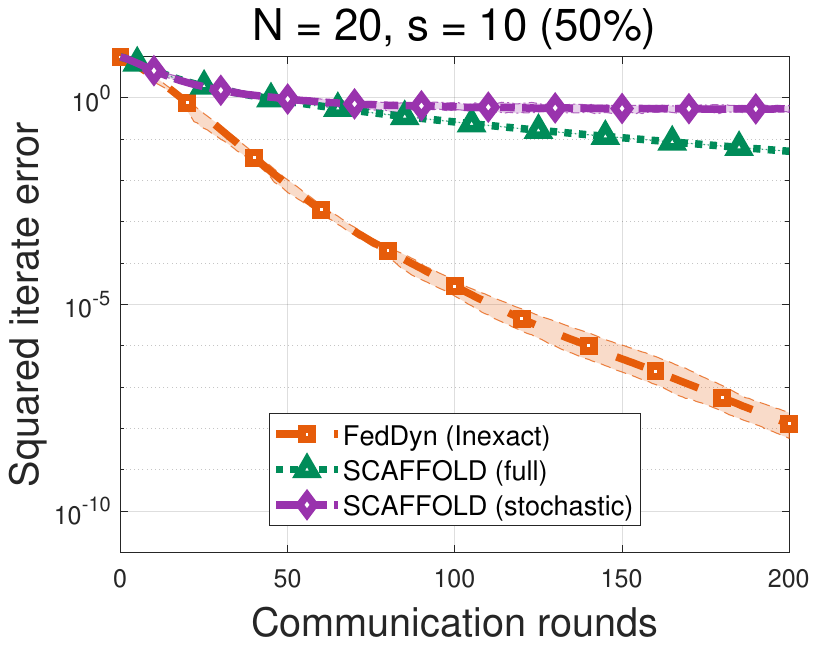}\hfill
    \includegraphics[trim={95pt 265pt 100pt 270pt},clip,width=0.32\linewidth]{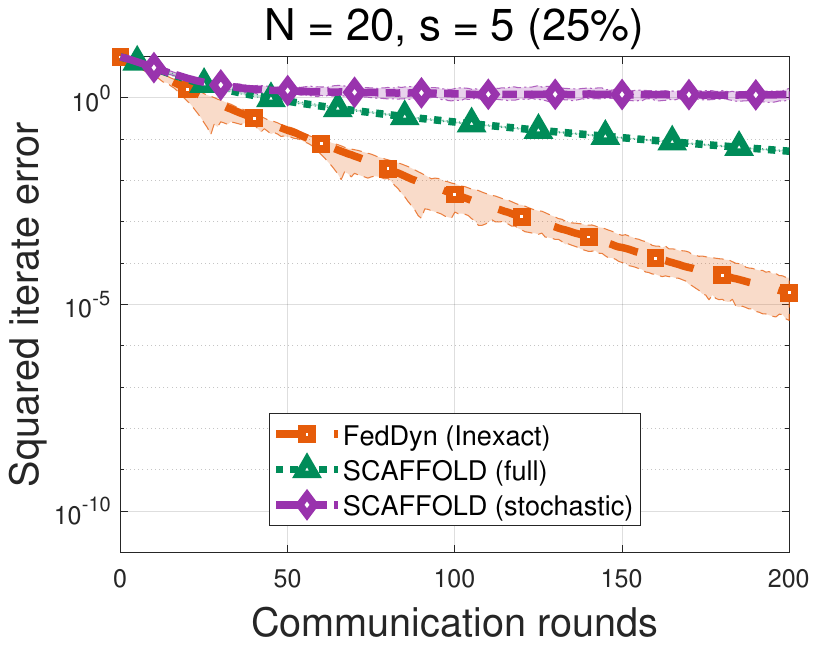}
    \par\medskip

    \includegraphics[trim={95pt 265pt 100pt 270pt},clip,width=0.32\linewidth]{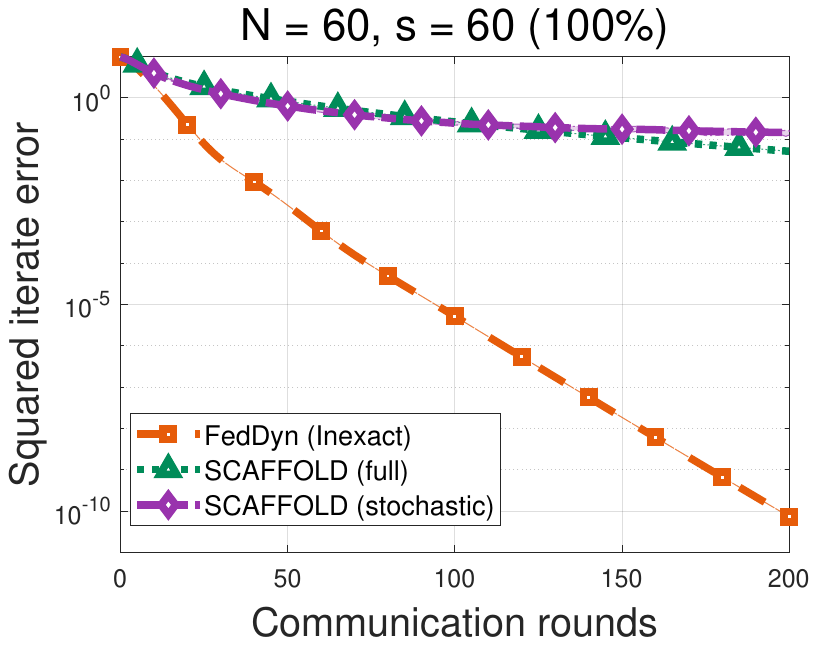}\hfill
    \includegraphics[trim={95pt 265pt 100pt 270pt},clip,width=0.32\linewidth]{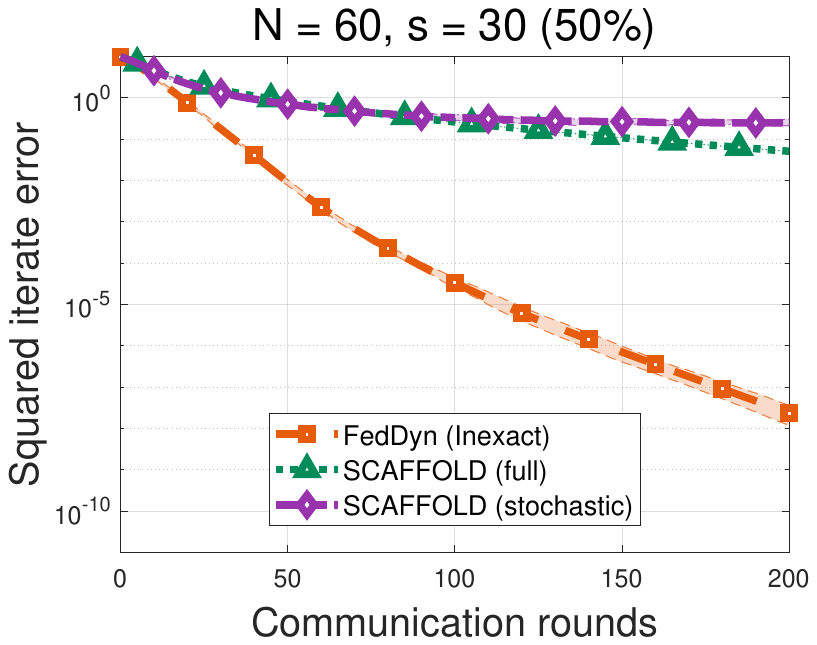}\hfill
    \includegraphics[trim={95pt 265pt 100pt 270pt},clip,width=0.32\linewidth]{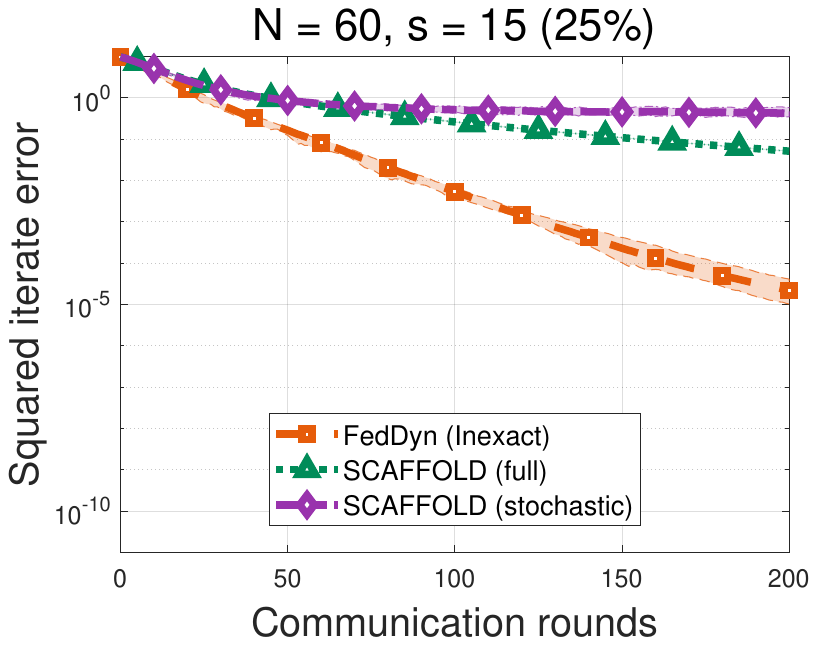}
    \par\medskip

    \includegraphics[trim={95pt 265pt 100pt 270pt},clip,width=0.32\linewidth]{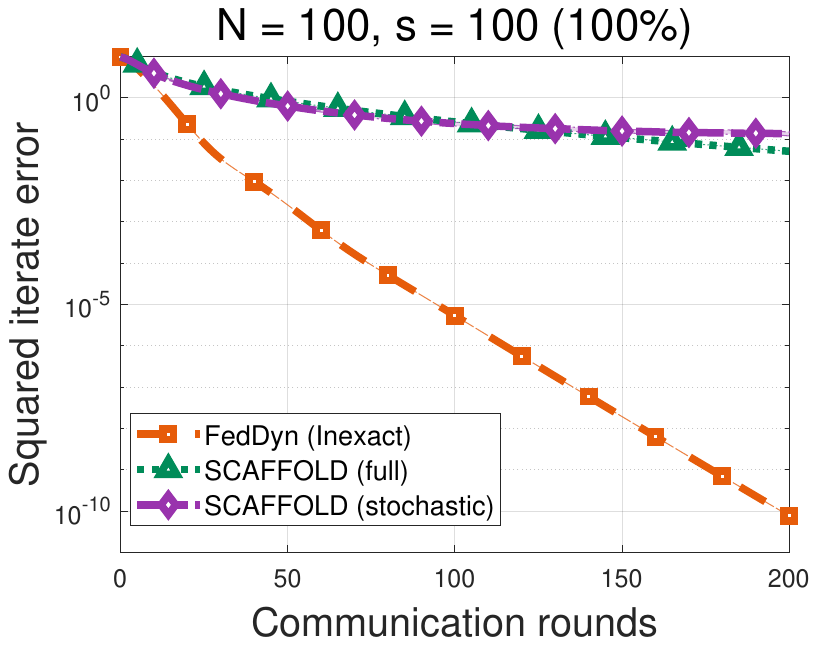}\hfill
    \includegraphics[trim={95pt 265pt 100pt 270pt},clip,width=0.32\linewidth]{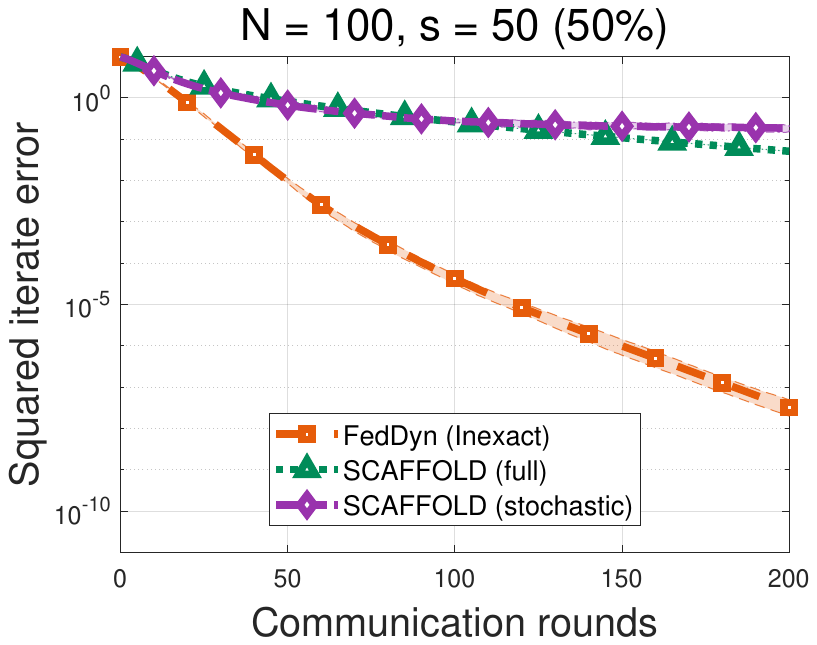}\hfill
    \includegraphics[trim={95pt 265pt 100pt 270pt},clip,width=0.32\linewidth]{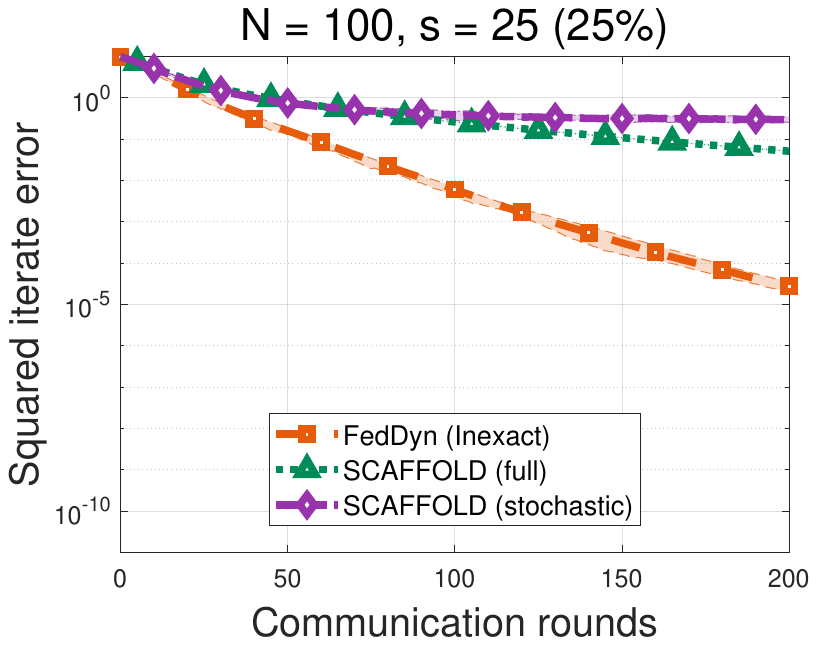}

    \caption{Squared iterate error versus communication rounds for
    distributed logistic regression on MNIST digits 4 and 9.
    Rows correspond to \(N=20,60,100\) clients from top to bottom,
    and columns correspond to \(s=N,N/2,N/4\) participating clients
    per round from left to right. Curves show the mean over 20 trials,
    and shaded bands span the minimum and maximum errors at each round.}
    \label{fig:convergence_curve_2}
\end{figure}

We compute $x^\star$ by Newton's method with gradient-norm tolerance
$10^{-10}$. Figure~\ref{fig:convergence_curve_2} shows the mean squared
iterate error $\|x^r-x^\star\|^2$ over 20 trials, with shaded bands
spanning the minimum and maximum errors at each round.

Figure~\ref{fig:convergence_curve_2} shows that, for each $N$,
\fdynInext{} attains a smaller error at round~200 as the participation
rate increases. The error of \scafFull{} continues to decrease and is
similar across all nine configurations. For \scafSto{}, higher
participation reduces the error in the later rounds. Under full participation, increasing $N$ reduces the mean squared
iterate error of \scafSto{}, but the reduction at round~200 falls
short of $1/N$ scaling.


\section{Conclusions}
\label{sec:conclusion}

Under the averaged strong convexity assumption, we establish convergence guarantees for \fdyn (both exact and inexact subproblem solutions) and \scaf. Moreover, we show that \fdyn with exact subproblem solutions and \scaf with deterministic gradients achieve the same convergence complexity.
Several interesting directions remain open. First, it is unclear whether \fdyn achieves an improved dependence on the condition number, in a manner analogous to the guarantees in~\cite{condat2024simple}. Second, designing new inexact variants of \fdyn that preserve the favorable convergence behavior of the exact method remains an important practical and theoretical challenge. 

\backmatter

\bmhead{Acknowledgements}

This work was partially supported by the Guangdong Major Project of Basic Research [2026B0303000015], the National Natural Science Foundation of China [72495131, 82441027], Guangdong Provincial Key Laboratory of Mathematical Foundations for Artificial Intelligence [2023B1212010001, ZSZDSYS20251011093302004].

\begin{data_availability}
    The MNIST dataset used in Section~\ref{sec:logisticExperiment} to classify handwritten digits 4 and 9 is available at \url{https://github.com/cvdfoundation/mnist}.
\end{data_availability}

\section*{Declarations}

\begin{conflicts_of_interest}
    No potential conflict of interest was reported by the authors.
\end{conflicts_of_interest}


\begin{appendices}

\section{Preliminary Lemmas}
\label{sec:preliminary-lemmas}
\begin{proof}[Proof of Lemma ~\ref{lem:mainineq}]
Let $f$ be convex and $L$-smooth. By \(L\)-smoothness~\cite{nesterov2004introductory}, we have
\begin{equation*}
    \langle \nabla f(z), x - z \rangle \geq f(x)- f(z) - \frac{L}{2} \| x - z \|^2 
\end{equation*}
for any $x, z$ and
\begin{equation*}
    \langle \nabla f(z),z-y\rangle
\geq
f(z)-f(y)+\frac{1}{2L}\|\nabla f(z)-\nabla f(y)\|^2
\end{equation*}
for any $y, z$. Adding these two inequalities gives inequality~\eqref{eq:mainineq}.
\end{proof}

\begin{lemma}[Young's Inequality with a Parameter]
    \label{lem:young}
    For any vectors $a, b \in \mathbb{R}^d$ and any constant $\alpha > 0$, we have
    \begin{equation*}
        \langle a, b \rangle \leq \frac{1}{2\alpha} \| a \|^2 + \frac{\alpha}{2} \| b \|^2.
    \end{equation*}
    The commonly used form $\langle a, b \rangle \leq \frac{1}{2} \| a \|^2 + \frac{1}{2} \| b \|^2$ is recovered by setting $\alpha = 1$.
\end{lemma}

\section{Convergence of \fdyn under \asc}
\label{app:feddyn}

{\bf Notation for Optimality Deviations.} 
For convenience, we measure all iterates relative to the optimal point $x^*$, the (unique) minimizer of $\bar f$, and define
\[
\Delta x_i^{r} \coloneqq x_i^{r}-x^*,\qquad
\Delta \lambda_i^{r} \coloneqq \lambda_i^{r}-\nabla f_i(x^*),\qquad
\Delta h^{r} \coloneqq h^{r}-\nabla \bar f(x^*).
\]
Since $\nabla \bar f(x^*)=0$, we have $\Delta h^{r}=h^{r}$. For the aggregated iterate
\[
x^{r} \coloneqq \frac{1}{s}\sum_{i\in S^{r-1}} x_i^{r},
\qquad
\Delta x^{r} \coloneqq x^{r}-x^*.
\]
Moreover, by the definition of $h^r$,
\[
\frac{1}{N}\sum_{i=1}^N \Delta \lambda_i^{r} = \Delta h^{r}.
\]

We first record a basic sampling bound that will be used later.
\begin{lemma}
\label{lem:feddyn_boundvariance}
Let $S^r\subseteq [N]$ be sampled uniformly at random among all subsets of cardinality $s$. Then
    \begin{equation}
        \label{eq:variance_bound_main}
            \mathbb{E}\Bigg[
            \sum_{i\in S^r}
            \|\lambda_i^r - h^r -\lambda_i^{r+1}\|^2
            \,\Big|\, \mathcal{F}_r
            \Bigg]
            \leq 
            \frac{2s}{N}
            \sum_{i=1}^N \|\Delta\lambda_i^r \|^2 + 2\,
            \mathbb{E}\Bigg[
            \sum_{i\in S^r}
            \|\Delta\lambda_i^{r+1} \|^2
            \,\Big|\, \mathcal{F}_r
            \Bigg].
    \end{equation}
\end{lemma}

\begin{proof}
Using $\Delta\lambda_i^r=\lambda_i^r-\nabla f_i(x^*)$ and $\Delta h^r=h^r-\nabla \bar f(x^*)=h^r$, we have
\[
\lambda_i^r-h^r-\lambda_i^{r+1}
=\Delta\lambda_i^r-\Delta h^r-\Delta\lambda_i^{r+1}.
\]
By Young's inequality, it follows that
\begin{equation}\label{eq:pointwise_bound_clean}
\|\lambda_i^r-h^r-\lambda_i^{r+1}\|^2
\le
2\|\Delta\lambda_i^r-\Delta h^r\|^2
+2\|\Delta\lambda_i^{r+1}\|^2.
\end{equation}
Expanding the first term gives
\[
\|\Delta\lambda_i^r-\Delta h^r\|^2
=
\|\Delta\lambda_i^r\|^2
-2\langle \Delta\lambda_i^r,\Delta h^r\rangle
+\|\Delta h^r\|^2.
\]

Summing \eqref{eq:pointwise_bound_clean} over $i\in S^r$ and taking conditional expectation given $\mathcal{F}_r$, we use that $S^r$ is sampled
uniformly and independently of $\mathcal{F}_r$, so
\begin{equation*}
\mathbb{E}\!\left[\sum_{i\in S^r}\|\Delta\lambda_i^r\|^2\mid \mathcal{F}_r\right]
=\frac{s}{N}\sum_{i=1}^N\|\Delta\lambda_i^r\|^2,
\quad
\mathbb{E}\!\left[\sum_{i\in S^r}\Delta\lambda_i^r\mid \mathcal{F}_r\right]
=\frac{s}{N}\sum_{i=1}^N\Delta\lambda_i^r=s\Delta h^r.
\end{equation*}

Plugging these identities into the conditional expectation of the sum of \eqref{eq:pointwise_bound_clean} yields
\begin{align*}
&\mathbb{E}\Bigg[
\sum_{i\in S^r}\|\lambda_i^r-h^r-\lambda_i^{r+1}\|^2
\;\Big|\;\mathcal{F}_r
\Bigg]\\
\le&
2\left[
\frac{s}{N}\sum_{i=1}^N\|\Delta\lambda_i^r\|^2
-2s\|\Delta h^r\|^2
+s\|\Delta h^r\|^2
\right]
+2\,\mathbb{E}\Bigg[
\sum_{i\in S^r}\|\Delta\lambda_i^{r+1}\|^2
\;\Big|\;\mathcal{F}_r
\Bigg] \\
\le&
\frac{2s}{N}\sum_{i=1}^N\|\Delta\lambda_i^r\|^2
+2\,\mathbb{E}\Bigg[
\sum_{i\in S^r}\|\Delta\lambda_i^{r+1}\|^2
\;\Big|\;\mathcal{F}_r
\Bigg],
\end{align*}
which is exactly \eqref{eq:variance_bound_main}.

\end{proof}

\begin{proof}[Proof of Theorem~\ref{thm:fdyn}]
We prove a one-step contraction of a suitable Lyapunov function. Throughout, condition on $\mathcal F_r$ and use the deviation notation
$\Delta x^r=x^r-x^*$ and $\Delta\lambda_i^r=\lambda_i^r-\nabla f_i(x^*)$.

For each sampled client $i\in S^r$, the update satisfies
\[
x_i^{r+1}=y^r+\frac{1}{\gamma}(\lambda_i^r-\lambda_i^{r+1}),
\qquad
y^r=x^r-\frac{1}{\gamma}h^r,
\]
hence
\[
\Delta x_i^{r+1}
=\Delta x^r+\frac{1}{\gamma}(\lambda_i^r-h^r-\lambda_i^{r+1}).
\]
Expanding the square gives
\[
\|\Delta x_i^{r+1}\|^2
=
\|\Delta x^r\|^2
+\frac{2}{\gamma}\langle \lambda_i^r-h^r-\lambda_i^{r+1},\Delta x^r\rangle
+\frac{1}{\gamma^2}\|\lambda_i^r-h^r-\lambda_i^{r+1}\|^2.
\]
Summing over $i\in S^r$ and taking $\mathbb E[\cdot\mid \mathcal F_r]$, the term involving $\lambda_i^r-h^r$ vanishes because
$\mathbb E\!\left[\sum_{i\in S^r}(\lambda_i^r-h^r)\mid\mathcal F_r\right]=0$ (uniform sampling and $h^r=\frac1N\sum_i\lambda_i^r$), yielding
\begin{equation}\label{eq:sum_expand_clean}
\begin{aligned}
\mathbb{E}\Bigl[ \sum_{i \in S^r} \|\Delta x_i^{r+1}\|^2 \,\Bigm|\, \mathcal{F}_r \Bigr]
=& s \|\Delta x^r\|^2
-\frac{2}{\gamma}\mathbb{E}\Bigl[ \sum_{i \in S^r} \langle \lambda_i^{r+1}, \Delta x^r \rangle \,\Bigm|\, \mathcal{F}_r \Bigr] \\
& + \frac{1}{\gamma^2}\mathbb{E}\Bigl[ \sum_{i \in S^r} \|\lambda_i^r - h^r - \lambda_i^{r+1} \|^2 \,\Bigm|\, \mathcal{F}_r \Bigr].
\end{aligned}
\end{equation}

\textbf{Control of the inner product term.} Let $z_i^{r+1}$ denote the \emph{exact} minimizer of the local subproblem (as defined earlier), and define the subproblem error
$e_i^{r+1}\coloneqq x_i^{r+1}-z_i^{r+1}$. Then the optimality condition of the subproblem gives
\[
\nabla f_i(z_i^{r+1})=\lambda_i^r-\gamma(z_i^{r+1}-y^r),
\qquad\Rightarrow\qquad
\lambda_i^{r+1}=\nabla f_i(z_i^{r+1})-\gamma e_i^{r+1}.
\]
Apply Lemma~\ref{lem:mainineq} to $f_i$ with $(x,y,z)=(x^r,x^*,z_i^{r+1})$:
\[
\langle \nabla f_i(z_i^{r+1}),\Delta x^r\rangle
\ge
f_i(x^r)-f_i(x^*)
+\frac{1}{2L}\|\nabla f_i(z_i^{r+1})-\nabla f_i(x^*)\|^2
-\frac{L}{2}\|x^r-z_i^{r+1}\|^2 .
\]
Multiplying by $-1$, substituting $\lambda_i^{r+1}=\nabla f_i(z_i^{r+1})-\gamma e_i^{r+1}$, and summing over $i\in S^r$ yields
\begin{equation}\label{eq:innerprod_bound_clean}
\begin{aligned}
-\sum_{i\in S^r}\langle \lambda_i^{r+1},\Delta x^r\rangle
\le{}&
-\sum_{i\in S^r}\bigl(f_i(x^r)-f_i(x^*)\bigr)
-\frac{1}{2L}\sum_{i\in S^r}\|\nabla f_i(z_i^{r+1})-\nabla f_i(x^*)\|^2\\
&+\frac{L}{2}\sum_{i\in S^r}\|x^r-z_i^{r+1}\|^2
+\gamma\sum_{i\in S^r}\langle e_i^{r+1},\Delta x^r\rangle .
\end{aligned}
\end{equation}
Taking conditional expectation and using uniform sampling gives
\[
\mathbb{E}\Bigl[\sum_{i\in S^r}(f_i(x^r)-f_i(x^*))\mid\mathcal F_r\Bigr]
=s\bigl(\bar f(x^r)-\bar f(x^*)\bigr)
\ge \frac{s\mu}{2}\|\Delta x^r\|^2,
\]
where the last inequality uses $\mu$-strong convexity of $\bar f$.

Using
\[
z_i^{r+1}-x^r=\frac{1}{\gamma}(\lambda_i^r-h^r-\lambda_i^{r+1})-e_i^{r+1},
\qquad
\nabla f_i(z_i^{r+1})-\nabla f_i(x^*)
=\Delta\lambda_i^{r+1}+\gamma e_i^{r+1},
\]
and Young's inequality with parameters $\theta_1,\theta_2,\theta_3>0$, we obtain
\begin{align}
\|x^r-z_i^{r+1}\|^2
&\le (1+\theta_1)\frac{1}{\gamma^2}\|\lambda_i^r-h^r-\lambda_i^{r+1}\|^2
+\Bigl(1+\frac{1}{\theta_1}\Bigr)\|e_i^{r+1}\|^2, \label{eq:bound_xz}\\
-\|\nabla f_i(z_i^{r+1})-\nabla f_i(x^*)\|^2
&\le -\frac{1}{1+\theta_2}\|\Delta\lambda_i^{r+1}\|^2+\frac{\gamma^2}{\theta_2}\|e_i^{r+1}\|^2, \label{eq:bound_grad}\\
2\langle e_i^{r+1},\Delta x^r\rangle
&\le \theta_3\|\Delta x^r\|^2+\frac{1}{\theta_3}\|e_i^{r+1}\|^2. \label{eq:bound_cross}
\end{align}
(When subproblems are solved exactly, $e_i^{r+1}=0$ and these error terms disappear.)

Substitute \eqref{eq:innerprod_bound_clean}--\eqref{eq:bound_cross} into \eqref{eq:sum_expand_clean}, and impose $L/\gamma\le 1/2$ to simplify constants.
This yields
\begin{equation*}
\begin{aligned}
&\mathbb{E}\Bigl[ \sum_{i \in S^r} \|\Delta x_i^{r+1}\|^2 \,\Bigm|\, \mathcal{F}_r \Bigr]
+ \frac{1}{(1+\theta_2)\gamma L}\,
\mathbb{E}\Bigl[ \sum_{i \in S^r} \|\Delta\lambda_i^{r+1}\|^2 \,\Bigm|\, \mathcal{F}_r \Bigr] \\
\le{}&
\Bigl(1-\frac{\mu}{\gamma}+\theta_3\Bigr)s\|\Delta x^r\|^2
+\frac{3+\theta_1}{2\gamma^2}\,
\mathbb{E}\Bigl[ \sum_{i \in S^r} \|\lambda_i^r-h^r-\lambda_i^{r+1}\|^2 \,\Bigm|\, \mathcal{F}_r \Bigr] \\
&+\Bigl[\frac12\Bigl(1+\frac{1}{\theta_1}\Bigr)+\frac{\gamma}{\theta_2L}+\frac{1}{\theta_3}\Bigr]\,
\mathbb{E}\Bigl[\sum_{i \in S^r}\|e_i^{r+1}\|^2 \,\Bigm|\, \mathcal{F}_r \Bigr].
\end{aligned}
\end{equation*}
By Lemma~\ref{lem:feddyn_boundvariance},
\[
\mathbb{E}\Bigl[ \sum_{i \in S^r} \|\lambda_i^r-h^r-\lambda_i^{r+1}\|^2 \,\Bigm|\, \mathcal{F}_r \Bigr]
\le
\frac{2s}{N}\sum_{i=1}^N\|\Delta\lambda_i^r\|^2
+2\mathbb{E}\Bigl[\sum_{i\in S^r}\|\Delta\lambda_i^{r+1}\|^2 \,\Bigm|\, \mathcal{F}_r\Bigr].
\]
Also, since $\lambda_i^{r+1}=\lambda_i^r$ for $i\notin S^r$,
\[
\mathbb{E}\Bigl[\sum_{i\in S^r}\|\Delta\lambda_i^{r+1}\|^2 \,\Bigm|\, \mathcal{F}_r\Bigr]
=
\mathbb{E}\Bigl[\sum_{i=1}^N\|\Delta\lambda_i^{r+1}\|^2 \,\Bigm|\, \mathcal{F}_r\Bigr]
-\frac{N-s}{N}\sum_{i=1}^N\|\Delta\lambda_i^{r}\|^2.
\]
By Jensen's inequality,
\[
\|\Delta x^{r+1}\|^2=\left\|\frac{1}{s}\sum_{i\in S^r}\Delta x_i^{r+1}\right\|^2
\le \frac{1}{s}\sum_{i\in S^r}\|\Delta x_i^{r+1}\|^2,
\]
hence
\[
s\,\mathbb{E}\bigl[\|\Delta x^{r+1}\|^2\mid\mathcal F_r\bigr]
\le
\mathbb{E}\Bigl[\sum_{i\in S^r}\|\Delta x_i^{r+1}\|^2\mid\mathcal F_r\Bigr].
\]

Collecting the previous bounds and defining
\[
A_1 \coloneqq \frac{1}{s}\Bigl[\frac{1}{(1+\theta_2)\gamma L}-\frac{3+\theta_1}{\gamma^2}\Bigr],
\qquad
\psi^r \coloneqq \|\Delta x^r\|^2 + A_1\sum_{i=1}^N\|\Delta\lambda_i^r\|^2,
\]
we obtain the recursion
\[
\mathbb{E}[\psi^{r+1}\mid\mathcal F_r]
\le
\Bigl(1-\frac{\mu}{\gamma}+\theta_3\Bigr)\psi^r
+
\frac{1}{s}\Bigl[\frac12\Bigl(1+\frac{1}{\theta_1}\Bigr)+\frac{\gamma}{\theta_2L}+\frac{1}{\theta_3}\Bigr]
\mathbb{E}\Bigl[\sum_{i\in S^r}\|e_i^{r+1}\|^2\mid\mathcal F_r\Bigr],
\]
where 
\begin{equation}\label{eq:lambda_coeff_condition}
\Bigl(1-\frac{\mu}{\gamma}+\theta_3\Bigr)A_1
\ge
\frac{3+\theta_1}{\gamma^2}\frac{1}{N}+A_1\frac{N-s}{N}.
\end{equation}
is applied.

(i) For the inexact-solve case, we set $\theta_3=\mu/(2\gamma)$ and use the assumed error model
$\mathbb{E}\bigl[\sum_{i\in S^r}\|e_i^{r+1}\|^2\mid\mathcal F_r\bigr]\le s(\epsilon^{r+1})^2$
to obtain
\[
\mathbb{E}[\psi^{r+1}\mid\mathcal F_r]
\le
\Bigl(1-\frac{\mu}{2\gamma}\Bigr)\psi^r
+
A_2(\epsilon^{r+1})^2,
\qquad
A_2=\frac12\Bigl(1+\frac{1}{\theta_1}\Bigr)+\frac{\gamma}{\theta_2L}+\frac{2\gamma}{\mu}.
\]
The stepsize conditions in Theorem~\ref{thm:fdyn}(i) ensure $A_1>0$ and that~\eqref{eq:lambda_coeff_condition} is valid.

(ii) For the exact-solve case, $e_i^{r+1}\equiv 0$, and we may take $\theta_1=\theta_2=\theta_3=0$, which yields
\[
\mathbb{E}[\psi^{r+1}\mid\mathcal F_r]\le \Bigl(1-\frac{\mu}{\gamma}\Bigr)\psi^r,
\qquad
A_1=\frac{1}{s}\Bigl(\frac{1}{\gamma L}-\frac{3}{\gamma^2}\Bigr),
\]
and the stepsize condition in Theorem~\ref{thm:fdyn}(ii) guarantees $A_1>0$, \eqref{eq:lambda_coeff_condition} is valid, and $1-\mu/\gamma\in(0,1)$.
This completes the proof. 

\end{proof}


\section{Convergence of \scaf under \asc}
\label{app:scaffold}

We first present several auxiliary results that will be used in the analysis.

\begin{lemma}[Lemma 4 in~\cite{karimireddy2020scaffold}]\label{lem:scaffold_random_seperate}
Let $\{\Xi_i\}_{i=1}^\tau\subset \mathbb{R}^d$ be (not necessarily independent) random vectors, and let $\{\xi_i\}_{i=1}^\tau\subset \mathbb{R}^d$ satisfy
\[
\mathbb{E}[\Xi_1]=\xi_1,
\qquad
\mathbb{E}[\Xi_i\mid \Xi_1,\ldots,\Xi_{i-1}]=\xi_i,\quad i=2,\ldots,\tau,
\]
so that $\{\Xi_i-\xi_i\}$ is a martingale difference sequence. Assume moreover that
\[
\mathbb{E}\!\left[\|\Xi_1-\xi_1\|^2\right]\le \sigma^2,
\qquad
\mathbb{E}\!\left[\|\Xi_i-\xi_i\|^2\mid \Xi_1,\ldots,\Xi_{i-1}\right]\le \sigma^2,\quad i=2,\ldots,\tau.
\]
Then,
\[
\mathbb{E}\Bigl[\Bigl\|\sum_{i=1}^\tau \Xi_i\Bigr\|^2\Bigr]
\le
2\,\mathbb{E}\Bigl[\Bigl\|\sum_{i=1}^\tau \xi_i\Bigr\|^2\Bigr]
+ 2\tau\sigma^2.
\]
\end{lemma}

\begin{proof}
Write $\Xi_i=\xi_i+(\Xi_i-\xi_i)$ and apply $\|a+b\|^2\le 2\|a\|^2+2\|b\|^2$:
\begin{equation}\label{eq:mds_split}
\mathbb{E}\Bigl[\Bigl\|\sum_{i=1}^\tau \Xi_i\Bigr\|^2\Bigr]
\le
2\,\mathbb{E}\Bigl[\Bigl\|\sum_{i=1}^\tau \xi_i\Bigr\|^2\Bigr]
+
2\,\mathbb{E}\Bigl[\Bigl\|\sum_{i=1}^\tau (\Xi_i-\xi_i)\Bigr\|^2\Bigr].
\end{equation}
For the second term, expand the square:
\[
\mathbb{E}\Bigl[\Bigl\|\sum_{i=1}^\tau (\Xi_i-\xi_i)\Bigr\|^2\Bigr]
=
\sum_{i=1}^\tau \mathbb{E}\|\Xi_i-\xi_i\|^2
+
2\sum_{1\le i<j\le \tau}\mathbb{E}\langle \Xi_i-\xi_i,\Xi_j-\xi_j\rangle.
\]
The cross terms vanish: for $i<j$, conditioning on $(\Xi_1,\ldots,\Xi_{j-1})$ gives
\[
\mathbb{E}\bigl[\langle \Xi_i-\xi_i,\Xi_j-\xi_j\rangle \mid \Xi_1,\ldots,\Xi_{j-1}\bigr]
=
\bigl\langle \Xi_i-\xi_i,\mathbb{E}[\Xi_j-\xi_j\mid \Xi_1,\ldots,\Xi_{j-1}]\bigr\rangle
=0,
\]
hence $\mathbb{E}\langle \Xi_i-\xi_i,\Xi_j-\xi_j\rangle=0$. Therefore,
\[
\mathbb{E}\Bigl[\Bigl\|\sum_{i=1}^\tau (\Xi_i-\xi_i)\Bigr\|^2\Bigr]
=
\sum_{i=1}^\tau \mathbb{E}\|\Xi_i-\xi_i\|^2
\le
\sum_{i=1}^\tau \sigma^2
=
\tau\sigma^2,
\]
where we used the bounded (conditional) second-moment assumption and the tower rule.
Substituting into \eqref{eq:mds_split} proves the claim. 

\end{proof}

\smallskip
\noindent\textbf{Notation for Optimality Deviations.}
To simplify the presentation, we measure all quantities relative to the optimizer $y^*$ and define
\begin{equation*}
    \Delta y^r \coloneqq y^r - y^*, \qquad
    \Delta \lambda_i^r \coloneqq \lambda_i^r - \nabla f_i(y^*), \qquad
    \Delta h^r \coloneqq h^r - \nabla \bar f(y^*).
\end{equation*}
Note that $\Delta y^r$ denotes the deviation of the global iterate from $y^*$; it should not be confused with the algorithmic quantity $[\Delta y]^r$
defined within \scaf.

\begin{lemma}
    \label{lem:scaffold_bound}
    Under Assumption~\ref{ass:BoundStoGrad}, for every communication round $r$,    
    \begin{align*}
        \mathbb{E}\bigl[\|\eta_g[\Delta y]^r\|^2 \mid \tilde{\mathcal{F}}_r\bigr]
        \leq
        K^2\eta_l^2\eta_g^2\biggl[
        \frac{2}{N}\sum_{i=1}^N\|\Delta \lambda_i^r\|^2
        +2\,\mathbb{E}\bigl[\frac{1}{s}\sum_{i\in S^r}\|\Delta \lambda_i^{r+1}\|^2 \mid \tilde{\mathcal{F}}_r\bigr]
        \biggr].
    \end{align*}
\end{lemma}

\begin{proof}
From Algorithm~\ref{alg:Scaffold}, for $i\in S^r$ we have
\[
\lambda_i^{r+1}
=\lambda_i^r-h^r+\frac{1}{K\eta_l}(y^r-x_i^{r,K})
=\frac{1}{K}\sum_{k=0}^{K-1} g_i(x_i^{r,k}),
\]
and therefore
\[
[\Delta x]_i^r
=x_i^{r,K}-y^r
=-K\eta_l\Bigl(\frac{1}{K}\sum_{k=0}^{K-1} g_i(x_i^{r,k})-\lambda_i^r+h^r\Bigr).
\]
Averaging over $i\in S^r$ and multiplying by $\eta_g$ gives
\begin{equation}\label{eq:scaf_deltay_repr}
\eta_g[\Delta y]^r
=
-K\eta_l\eta_g\cdot
\frac{1}{s}\sum_{i\in S^r}
\Bigl(\frac{1}{K}\sum_{k=0}^{K-1} g_i(x_i^{r,k})-\lambda_i^r+h^r\Bigr).
\end{equation}
Since $y^*$ minimizes $\bar f$, we have $\nabla\bar f(y^*)=0$. Define, for $i\in S^r$,
\[
u_i
\coloneqq
\frac{1}{K}\sum_{k=0}^{K-1} g_i(x_i^{r,k})-\lambda_i^r+h^r-\nabla\bar f(y^*).
\]
Then \eqref{eq:scaf_deltay_repr} and Jensen's inequality imply
\begin{align}
\|\eta_g[\Delta y]^r\|^2
=K^2\eta_l^2\eta_g^2\Bigl\|\frac{1}{s}\sum_{i\in S^r}u_i\Bigr\|^2
\le
K^2\eta_l^2\eta_g^2\cdot \frac{1}{s}\sum_{i\in S^r}\|u_i\|^2. \label{eq:bound_1}
\end{align}

Next, decompose $u_i$ as
\[
u_i
=
\underbrace{\Bigl(\frac{1}{K}\sum_{k=0}^{K-1} g_i(x_i^{r,k})-\nabla f_i(y^*)\Bigr)}_{\eqqcolon\,a_i}
+
\underbrace{\bigl(\nabla f_i(y^*)-\lambda_i^r+h^r-\nabla\bar f(y^*)\bigr)}_{=-(\Delta\lambda_i^r-\Delta h^r)\eqqcolon\,b_i}.
\]
By Young's inequality,
\begin{equation}\label{eq:ui_split}
\|u_i\|^2
\le
2\|\Delta\lambda_i^r-\Delta h^r\|^2
+2\Bigl\|\frac{1}{K}\sum_{k=0}^{K-1} g_i(x_i^{r,k})-\nabla f_i(y^*)\Bigr\|^2.
\end{equation}
Taking conditional expectation and using uniform sampling gives
\[
\mathbb{E}\Bigl[\frac{1}{s}\sum_{i\in S^r}\|\Delta\lambda_i^r-\Delta h^r\|^2 \,\Big|\, \tilde{\mathcal{F}}_r\Bigr]
=
\frac{1}{N}\sum_{i=1}^N\|\Delta\lambda_i^r-\Delta h^r\|^2
\le\;
\frac{1}{N}\sum_{i=1}^N\|\Delta\lambda_i^r\|^2,
\]
where we used $\Delta h^r=\frac{1}{N}\sum_{i=1}^N\Delta\lambda_i^r$. 
Moreover, since $\lambda_i^{r+1}=\frac{1}{K}\sum_{k=0}^{K-1} g_i(x_i^{r,k})$ for $i\in S^r$,
\[
\mathbb{E}\Bigl[\frac{1}{s}\sum_{i\in S^r}\Bigl\|\frac{1}{K}\sum_{k=0}^{K-1} g_i(x_i^{r,k})-\nabla f_i(y^*)\Bigr\|^2 \,\Big|\, \tilde{\mathcal{F}}_r\Bigr]
=
\mathbb{E}\Bigl[\frac{1}{s}\sum_{i\in S^r}\|\Delta\lambda_i^{r+1}\|^2 \,\Big|\, \tilde{\mathcal{F}}_r\Bigr].
\]
Finally, taking conditional expectation in \eqref{eq:bound_1} and combining with \eqref{eq:ui_split} and the two bounds above yields the claim. 

\end{proof}

Because each client performs multiple local update steps per communication round, an additional ``drift'' term appears, reflecting the deviation of local iterates from the global reference point. We next derive a bound on this drift term.

\begin{lemma}
    \label{lem:scaffold_client_drift}
    Under Assumption~\ref{ass:BoundStoGrad}, suppose the local stepsize satisfies $0<\eta_l\le \frac{2}{L}$. Then, for every round $r$,
    \begin{align*}
        &\frac 1K\sum_{k=0}^{K-1}\mathbb{E}\bigl[\frac 1s\sum_{i\in S^r}\|x_i^{r,k}-y^r\|^2\big|\tilde{\mathcal{F}}_r\bigr] \\
        \leq&\,2(1+K)K\eta_l^2\frac 1N\sum_{i=1}^N\|\Delta\lambda_i^r\|^2+2(1+K)K\eta_l^2L^2\|\Delta y^r\|^2+K\eta_l^2\sigma^2.
    \end{align*}
\end{lemma}
\begin{proof}
Fix a round $r$ and define the local drift energy
\[
V_k \coloneqq \mathbb{E}\Bigl[\frac{1}{s}\sum_{i\in S^r}\|x_i^{r,k}-y^r\|^2 \,\Big|\,\tilde{\mathcal{F}}_r\Bigr],
\qquad k=0,1,\ldots,K,
\]
so that $V_0=0$ since $x_i^{r,0}=y^r$.

For $i\in S^r$, the local update is
\[
x_i^{r,k}=x_i^{r,k-1}-\eta_l\bigl(g_i(x_i^{r,k-1})-\lambda_i^r+h^r\bigr),\qquad k\ge 1.
\]
Let $\mathcal G_{i,k-1}\coloneqq \sigma(x_i^{r,0},\ldots,x_i^{r,k-1},S^r)\vee \tilde{\mathcal F}_r$. By Assumption~\ref{ass:BoundStoGrad},
\[
\mathbb{E}[g_i(x_i^{r,k-1})\mid \mathcal G_{i,k-1}] = \nabla f_i(x_i^{r,k-1}),
\quad
\mathbb{E}\bigl[\|g_i(x_i^{r,k-1})-\nabla f_i(x_i^{r,k-1})\|^2\mid \mathcal G_{i,k-1}\bigr]\le \sigma^2.
\]
Using the conditional variance decomposition, for each $i\in S^r$,
\begin{align*}
&\mathbb{E}\bigl[\|x_i^{r,k}-y^r\|^2 \mid \mathcal G_{i,k-1}\bigr] \le
\|x_i^{r,k-1}-y^r-\eta_l(\nabla f_i(x_i^{r,k-1})-\lambda_i^r+h^r)\|^2+\eta_l^2\sigma^2 .
\end{align*}
Averaging over $i\in S^r$ and conditioning on $\tilde{\mathcal F}_r$ yields
\begin{equation}\label{eq:Vk_basic}
V_k
\le
\mathbb{E}\Bigl[\frac{1}{s}\sum_{i\in S^r}
\|x_i^{r,k-1}-y^r-\eta_l(\nabla f_i(x_i^{r,k-1})-\lambda_i^r+h^r)\|^2
\Bigm|\tilde{\mathcal F}_r\Bigr]
+\eta_l^2\sigma^2.
\end{equation}

Write
\[
\nabla f_i(x_i^{r,k-1})-\lambda_i^r+h^r
=
\bigl(\nabla f_i(x_i^{r,k-1})-\nabla f_i(y^r)\bigr)
+\bigl(\nabla f_i(y^r)-\lambda_i^r+h^r\bigr).
\]
Applying Young's inequality $\|a+b\|^2\le (1+\frac{1}{K})\|a\|^2+(1+K)\|b\|^2$ to \eqref{eq:Vk_basic} gives
\begin{equation}\label{eq:Vk_split}
V_k
\le
\Bigl(1+\frac{1}{K}\Bigr)\,T_k
+(1+K)\eta_l^2\,U
+\eta_l^2\sigma^2,
\end{equation}
where
\[
T_k \;\coloneqq\;
\mathbb{E}\Bigl[\frac{1}{s}\sum_{i\in S^r}
\|x_i^{r,k-1}-y^r-\eta_l(\nabla f_i(x_i^{r,k-1})-\nabla f_i(y^r))\|^2
\Bigm|\tilde{\mathcal F}_r\Bigr],
\]
and
\[
U \coloneqq
\frac{1}{N}\sum_{i=1}^N \|\lambda_i^r-\nabla f_i(y^r)-h^r\|^2
\quad
(\text{note that }U\text{ is }\tilde{\mathcal F}_r\text{-measurable}).
\]

For convex $L$-smooth $f_i$, the gradient mapping is nonexpansive for $0<\eta_l\le 2/L$, i.e.,
\[
\|u-\eta_l\nabla f_i(u) - (v-\eta_l\nabla f_i(v))\| \le \|u-v\|,\qquad \forall u,v.
\]
With $(u,v)=(x_i^{r,k-1},y^r)$ this implies $T_k \le V_{k-1}$. Hence, \eqref{eq:Vk_split} becomes
\begin{equation}\label{eq:Vk_recursion}
V_k
\le
\Bigl(1+\frac{1}{K}\Bigr)V_{k-1}
+(1+K)\eta_l^2\,U
+\eta_l^2\sigma^2.
\end{equation}

Iterating \eqref{eq:Vk_recursion} from $V_0=0$ yields, for $k=1,\ldots,K$,
\[
V_k
\le
\left(\sum_{t=0}^{k-1}\Bigl(1+\frac{1}{K}\Bigr)^t\right)
\Bigl((1+K)\eta_l^2\,U+\eta_l^2\sigma^2\Bigr).
\]
Averaging over $k=0,\ldots,K-1$ gives
\[
\frac{1}{K}\sum_{k=0}^{K-1}V_k
\le
\left(\frac{1}{K}\sum_{k=0}^{K-1}\sum_{t=0}^{k-1}\Bigl(1+\frac{1}{K}\Bigr)^t\right)
\Bigl((1+K)\eta_l^2\,U+\eta_l^2\sigma^2\Bigr).
\]
Moreover,
\[
\frac{1}{K}\sum_{k=0}^{K-1}\sum_{t=0}^{k-1}\Bigl(1+\frac{1}{K}\Bigr)^t
=
\sum_{k=0}^{K-1}\Bigl[\Bigl(1+\frac{1}{K}\Bigr)^k-1\Bigr]
=
K\Bigl(\Bigl(1+\tfrac{1}{K}\Bigr)^K-2\Bigr)
< K,
\]
where we used $\bigl(1+\tfrac{1}{K}\bigr)^K<e<3$. Therefore,
\begin{equation}\label{eq:Vk_avg_preU}
\frac{1}{K}\sum_{k=0}^{K-1}V_k
\le
K\Bigl((1+K)\eta_l^2\,U+\eta_l^2\sigma^2\Bigr).
\end{equation}

It remains to bound $U$. Using $\Delta\lambda_i^r=\lambda_i^r-\nabla f_i(y^*)$, $\Delta h^r=h^r-\nabla\bar f(y^*)$, and $\nabla\bar f(y^*)=0$, we have
\[
\lambda_i^r-\nabla f_i(y^r)-h^r
=
(\Delta\lambda_i^r-\Delta h^r) -\bigl(\nabla f_i(y^r)-\nabla f_i(y^*)\bigr).
\]
Thus, by $\|a+b\|^2\le 2\|a\|^2+2\|b\|^2$ and $L$-smoothness,
\[
U
\le
\frac{2}{N}\sum_{i=1}^N\|\Delta\lambda_i^r-\Delta h^r\|^2
+\frac{2}{N}\sum_{i=1}^N\|\nabla f_i(y^r)-\nabla f_i(y^*)\|^2
\le
\frac{2}{N}\sum_{i=1}^N\|\Delta\lambda_i^r\|^2
+2L^2\|\Delta y^r\|^2,
\]
where we used $\Delta h^r=\frac{1}{N}\sum_{i=1}^N\Delta\lambda_i^r$. Applying this bound into \eqref{eq:Vk_avg_preU} proves the lemma. 

\end{proof}

\begin{proof}[Proof of Theorem~\ref{thm:scaf}]
Fix a round $r$ and condition on $\tilde{\mathcal F}_r$.
Using $y^{r+1}=y^r+\eta_g[\Delta y]^r$, we have
\begin{equation}\label{eq:scaf_pf_expand}
\mathbb{E}\bigl[\|\Delta y^{r+1}\|^2\mid \tilde{\mathcal{F}}_r\bigr]
=
\|\Delta y^r\|^2
+2\eta_g\,\mathbb{E}\bigl[\langle[\Delta y]^r,\Delta y^r\rangle\mid \tilde{\mathcal{F}}_r\bigr]
+\mathbb{E}\bigl[\|\eta_g[\Delta y]^r\|^2\mid \tilde{\mathcal{F}}_r\bigr].
\end{equation}

We focus on the last two terms on the right-hand side.  From Algorithm~\ref{alg:Scaffold},
\[
[\Delta y]^r
=\frac{1}{s}\sum_{i\in S^r}(x_i^{r,K}-y^r)
=-\eta_l\frac{1}{s}\sum_{i\in S^r}\sum_{k=0}^{K-1}\bigl(g_i(x_i^{r,k})-\lambda_i^r+h^r\bigr).
\]
Taking conditional expectation and using unbiasedness of $g_i(\cdot)$, together with
$\mathbb{E}\bigl[\frac{1}{s}\sum_{i\in S^r}(\lambda_i^r-h^r)\mid \tilde{\mathcal F}_r\bigr]=0$, yields
\begin{equation*}
\mathbb{E}\bigl[[\Delta y]^r\mid \tilde{\mathcal{F}}_r\bigr]
=
-\eta_l\sum_{k=0}^{K-1}\mathbb{E}\biggl[\frac{1}{s}\sum_{i\in S^r}\nabla f_i(x_i^{r,k})\ \Bigm|\ \tilde{\mathcal{F}}_r\biggr].
\end{equation*}
Since each $f_i$ is convex and $L$-smooth, Lemma~\ref{lem:mainineq} with $(x,y,z)=(y^r,y^*,x_i^{r,k})$ implies
\[
-\langle \nabla f_i(x_i^{r,k}),\,\Delta y^r\rangle
\le
f_i(y^*)-f_i(y^r)
-\frac{1}{2L}\|\nabla f_i(x_i^{r,k})-\nabla f_i(y^*)\|^2
+\frac{L}{2}\|x_i^{r,k}-y^r\|^2.
\]
Substituting this bound into the inner-product term in \eqref{eq:scaf_pf_expand}, using uniform sampling, and applying Jensen's inequality to the gradient term yields (with $\tilde{\eta}\coloneqq K\eta_g\eta_l$)
\begin{equation*}
\begin{aligned}
2\eta_g\,\mathbb{E}\bigl[\langle[\Delta y]^r,\Delta y^r\rangle\mid \tilde{\mathcal{F}}_r\bigr]\le{}&
2\tilde{\eta}\bigl(\bar f(y^*)-\bar f(y^r)\bigr)-\frac{\tilde{\eta}}{L}\,G_r
 +\tilde{\eta}LD_r,
\end{aligned}
\end{equation*}
where
\begin{align*}
G_r &\coloneqq
\mathbb{E}\biggl[\frac{1}{s}\sum_{i\in S^r}\Bigl\|\frac{1}{K}\sum_{k=0}^{K-1}\nabla f_i(x_i^{r,k})-\nabla f_i(y^*)\Bigr\|^2 \Bigm| \tilde{\mathcal{F}}_r\biggr],\\
D_r &\coloneqq
\frac{1}{K}\sum_{k=0}^{K-1}\mathbb{E}\biggl[\frac{1}{s}\sum_{i\in S^r}\|x_i^{r,k}-y^r\|^2\Bigm|\tilde{\mathcal{F}}_r\biggr].
\end{align*}
Apply Lemma~\ref{lem:scaffold_random_seperate} conditionally (with $\tau=K$) to
$\frac{1}{K}\sum_{k=0}^{K-1}(g_i(x_i^{r,k})-\nabla f_i(y^*))=\Delta\lambda_i^{r+1}$ yields
\begin{equation*}
T_r\coloneqq\mathbb{E}\biggl[\frac{1}{s}\sum_{i\in S^r}\|\Delta\lambda_i^{r+1}\|^2 \,\Bigm|\, \tilde{\mathcal{F}}_r\biggr]
\le
2\,G_r
+\frac{2}{K}\sigma^2.
\end{equation*}
Thus we have 
\begin{equation}\label{eq:scaf_pf_after_innerprod}
\begin{aligned}
2\eta_g\,\mathbb{E}\bigl[\langle[\Delta y]^r,\Delta y^r\rangle\mid \tilde{\mathcal{F}}_r\bigr]\le{}&
2\tilde{\eta}\bigl(\bar f(y^*)-\bar f(y^r)\bigr)
-\frac 92\tilde{\eta}^2\,T_r \\
&+ (9\tilde{\eta}^2-\frac{\tilde{\eta}}{L})G_r + 9\tilde{\eta}^2\frac{\sigma^2}{K}+\tilde{\eta}L D_r.
\end{aligned}
\end{equation}
By choosing $\tilde{\eta}$ such that $9\tilde{\eta}^2-\frac{\tilde{\eta}}{L}\leq 0$, the coefficient of the term $G_r\geq 0$ is nonpositive, and the entire term can be dropped.

Moreover, Lemma~\ref{lem:scaffold_bound} implies
\begin{equation}\label{eq:scaf_pf_var}
\mathbb{E}\bigl[\|\eta_g[\Delta y]^r\|^2 \mid \tilde{\mathcal{F}}_r\bigr]
\le
\tilde{\eta}^2\biggl[
\frac{2}{N}\sum_{i=1}^N\|\Delta \lambda_i^r\|^2
+2\,T_r
\biggr].
\end{equation}

Combining \eqref{eq:scaf_pf_after_innerprod}--\eqref{eq:scaf_pf_var} with~\eqref{eq:scaf_pf_expand} and using $\bar f(y^r)-\bar f(y^*)\ge \frac{\mu}{2}\|\Delta y^r\|^2$, we obtain

\begin{equation}
\begin{split}
\mathbb{E}\bigl[\|\Delta y^{r+1}\|^2\mid\tilde{\mathcal{F}}_r\bigr]
\le{}&
\bigl(1-\tilde{\eta}\mu\bigr)\|\Delta y^r\|^2
-\frac 52\tilde{\eta}^2 T_r \\
&+\tilde{\eta}L\,D_r
+2\tilde{\eta}^2\cdot\frac{1}{N}\sum_{i=1}^N\|\Delta\lambda_i^r\|^2
+9\tilde{\eta}^2\frac{\sigma^2}{K}.
\end{split}\label{eq:tmp_before_add}
\end{equation}

Since $\lambda_i^{r+1}=\lambda_i^r$ for $i\notin S^r$, we have $\Delta\lambda_i^{r+1}=\Delta\lambda_i^r$ for $i\notin S^r$, and therefore
\[
T_r=\mathbb{E}\biggl[{1\over s}\sum_{i\in S^r}\|\Delta\lambda_i^{r+1}\|^2\ \Bigm|\ \tilde{\mathcal{F}}_r\biggr]
=
{1\over s}\sum_{i=1}^N\mathbb{E}\bigl[\|\Delta\lambda_i^{r+1}\|^2\mid\tilde{\mathcal{F}}_r\bigr]
-\frac{N-s}{Ns}\sum_{i=1}^N\|\Delta\lambda_i^r\|^2.
\]
Substituting into \eqref{eq:tmp_before_add} and defining
\[
B \coloneqq \frac 52\tilde{\eta}^2\frac {N}{s},
\qquad
\phi^r \coloneqq \|\Delta y^r\|^2 + \frac{B}{N}\sum_{i=1}^N\|\Delta\lambda_i^r\|^2,
\]
we obtain
\[
\mathbb{E}\bigl[\phi^{r+1}\mid\tilde{\mathcal{F}}_r\bigr]
\le \bigl(1-\tilde{\eta}\mu\bigr)\|\Delta y^r\|^2
 +\tilde{\eta}L\,D_r
+\Bigl(2\tilde{\eta}^2+\frac{N-s}{N}B\Bigr)\cdot\frac{1}{N}\sum_{i=1}^N\|\Delta\lambda_i^r\|^2
+9\tilde{\eta}^2\frac{\sigma^2}{K}.
\]

Then apply Lemma~\ref{lem:scaffold_client_drift} to bound the drift average $D_r$. This yields
\begin{equation*}
    \begin{split}
\mathbb{E}\bigl[\phi^{r+1}\mid\tilde{\mathcal{F}}_r\bigr]
\le&
\Bigl(1-\tilde{\eta}\bigl(\mu-2(1+K)K\eta_l^2L^3\bigr)\Bigr)\|\Delta y^r\|^2 +\tilde{\eta}\bigl(LK^2{\eta_l}^2+9\tilde{\eta}\bigr)\frac{\sigma^2}{K} \\
&+\Bigl(2\tilde{\eta}^2+\frac{N-s}{N}B+2\tilde{\eta}(1+K)K\eta_l^2L\Bigr)\frac{1}{N}\sum_{i=1}^N\|\Delta\lambda_i^r\|^2
.
\end{split}
\end{equation*}
Choose $\eta_l$ such that $\mu-2(1+K)K\eta_l^2L^3\ge \mu/2$,  so that the coefficient of $\|\Delta y^r\|^2 $ is bounded by $1-\tilde{\eta}\mu/2$.

To obtain a contraction in the Lyapunov function $\phi^r$, additionally require
\[
2\tilde{\eta}^2+\frac{N-s}{N}B+2\tilde{\eta}(1+K)K\eta_l^2L
\le
\Bigl(1-\frac{\tilde{\eta}\mu}{2}\Bigr)B,
\]
which is ensured by the second upper bound on $\tilde{\eta}$ assumed in the theorem.
Together with the preceding constraints on $\eta_l$ and $\tilde{\eta}$, the stepsize conditions stated in the theorem ensure the recursion
\[
\mathbb{E}\bigl[\phi^{r+1}\mid\tilde{\mathcal{F}}_r\bigr]
\le
\Bigl(1-\frac{\tilde{\eta}\mu}{2}\Bigr)\phi^r
+\tilde{\eta}\bigl(L\tilde{\eta}^2+9\tilde{\eta}\bigr)\frac{\sigma^2}{K}.
\]

This completes the proof.

\end{proof}

\section{Complexity of inexact FedDyn}
\label{app:fdyn-inexact-complexity}

We take \(0<\epsilon<1\) and use
\(\theta_1=\theta_2=1\), \(\gamma=16L+N\mu/s\), and
\(\kappa=L/\mu\). The inexact recursion gives
\begin{equation}
\label{eq:fdyn-unrolled}
\mathbb E[\psi^R]
\le \rho^R\psi^0+A_2\sum_{r=1}^R\rho^{R-r}(\epsilon^r)^2,
\qquad
\rho=1-\frac{\mu}{2\gamma},
\qquad
A_2=1+\frac{\gamma}{L}+\frac{2\gamma}{\mu}.
\end{equation}

\noindent\textbf{Local gradient evaluations.}
In round \(r\), client \(i\in S^{r-1}\) applies Nesterov's accelerated
gradient method to
\[
Q_i^{r-1}(u)
=f_i(u)-\langle\lambda_i^{r-1},u-y^{r-1}\rangle
+\frac{\gamma}{2}\|u-y^{r-1}\|^2.
\]
This objective is \(\gamma\)-strongly convex and
\((L+\gamma)\)-smooth, with \(\kappa_Q=(L+\gamma)/\gamma\le17/16\).
We initialize \(u_{-1}=u_0=y^{r-1}\) and use stepsize
\(1/(L+\gamma)\) and momentum
\((\sqrt{\kappa_Q}-1)/(\sqrt{\kappa_Q}+1)\).
If \(z_i^r\) denotes the exact minimizer, strong convexity gives
\[
\|u_k-z_i^r\|\le\frac{1}{\gamma}\|\nabla Q_i^{r-1}(u_k)\|.
\]
Thus, returning the first iterate satisfying
\(\|\nabla Q_i^{r-1}(u_k)\|\le\gamma\epsilon^r\)
ensures the required absolute accuracy.

Let \(g_i^r=\nabla Q_i^{r-1}(y^{r-1})
=\nabla f_i(y^{r-1})-\lambda_i^{r-1}\).
The standard accelerated gradient estimate~\cite{nesterov2004introductory},
together with smoothness and strong convexity, gives
\begin{equation}
\label{eq:fdyn-inner-gradient-rate}
\|\nabla Q_i^{r-1}(u_k)\|^2
\le 2\kappa_Q
\left(1-\frac{1}{\sqrt{\kappa_Q}}\right)^k\|g_i^r\|^2.
\end{equation}
Since \(\kappa_Q\le17/16\), the number \(K_i^r\) of local gradient
evaluations made by client \(i\) in round \(r\), including the
stopping tests, satisfies
\[
K_i^r=\mathcal O\!\left(
\log\!\left(2+\frac{\|g_i^r\|^2}
{\gamma^2(\epsilon^r)^2}\right)\right).
\]

To pass from a single client to the parallel count, observe that
\(y^{r-1}=x^{r-1}-h^{r-1}/\gamma\) and
\(h^{r-1}=N^{-1}\sum_i(\lambda_i^{r-1}-\nabla f_i(x^\star))\).
Lipschitz continuity and the triangle inequality yield
\[
\left(\frac1N\sum_i\|g_i^r\|^2\right)^{1/2}
\le L\|x^{r-1}-x^\star\|
+\left(1+\frac L\gamma\right)
\left(\frac1N\sum_i
\|\lambda_i^{r-1}-\nabla f_i(x^\star)\|^2\right)^{1/2}.
\]
By Cauchy--Schwarz and the definition of \(\psi^{r-1}\),
\[
\frac1N\sum_i\|g_i^r\|^2
\le
\left(L^2+\frac{(1+L/\gamma)^2}{NA_1}\right)\psi^{r-1}
=\mathcal O(L^2\psi^{r-1}),
\]
where the last bound follows from
\(1/(NA_1)=2(s/N)\gamma L/(1-8L/\gamma)\),
\(\gamma\ge16L\), and \(s\gamma/N\le17L\).
The maximum over the sampled clients is bounded by their sum.
Uniform sampling and Jensen's inequality therefore give
\begin{equation}
\label{eq:fdyn-parallel-round}
\mathbb E\!\left[\max_{i\in S^{r-1}}K_i^r\right]
=\mathcal O\!\left(
\log\!\left(2+
\frac{sL^2\,\mathbb E[\psi^{r-1}]}
     {\gamma^2(\epsilon^r)^2}\right)\right).
\end{equation}
We sum this bound over \(r\) to obtain the expected parallel gradient
count \(\mathbb E[G_{\rm par}(R)]\).

\begin{itemize}
    \item \textbf{Uniform inexactness.}
    Set \((\epsilon^r)^2=\mu\epsilon/(4\gamma A_2)\).
    Since \(1-\rho=\mu/(2\gamma)\), \eqref{eq:fdyn-unrolled} gives
    \[
    \mathbb E[\psi^R]\le \rho^R\psi^0+\frac{\epsilon}{2}.
    \]
    Using \(\rho^R\le\exp(-R\mu/(2\gamma))\) and
    \(\gamma/\mu=16\kappa+N/s\), we obtain the communication complexity
    \[
    R=\mathcal O\!\left(
    \left(\kappa+\frac Ns\right)\log\frac{1}{\epsilon}\right).
    \]
    Moreover, \(\mathbb E[\psi^{r-1}]\le \rho^{r-1}\psi^0+\epsilon/2\le \psi^0+\epsilon/2\) and
    \(A_2\le4\gamma/\mu\), so
    \(\gamma^2(\epsilon^r)^2\ge\mu^2\epsilon/16\).
    Substituting into \eqref{eq:fdyn-parallel-round} and summing gives
    \[
    \mathbb E[G_{\rm par}(R)]
    =\mathcal O\!\left(
    \left(\kappa+\frac Ns\right)\log\frac{1}{\epsilon}\,
    \log\!\left(2+\frac{s\kappa^2}{\epsilon}\right)\right).
    \]
    
    \item \textbf{Geometrically decaying inexactness.}
    Let \(\delta=\mu/(8\gamma)\) and set
    \((\epsilon^r)^2=(\delta/A_2)(\rho-\delta)^{r-1}\).
    Evaluating the geometric sum in \eqref{eq:fdyn-unrolled} yields
    \[
    \mathbb E[\psi^R]
    \le \rho^R\psi^0+
    \delta\sum_{r=1}^R\rho^{R-r}(\rho-\delta)^{r-1}
    =\rho^R\psi^0+\rho^R-(\rho-\delta)^R
    \le\rho^R(\psi^0+1).
    \]
    Thus the communication complexity is again
    \[
    R=\mathcal O\!\left(
    \left(\kappa+\frac Ns\right)\log\frac{1}{\epsilon}\right).
    \]
    For the local solves, \(A_2\le4\gamma/\mu\) gives
    \(\gamma^2(\epsilon^r)^2\ge
    (\mu^2/32)(\rho-\delta)^{r-1}\).
    Combining this with \(\mathbb E[\psi^{r-1}]\le (1+\psi^0)\rho^{r-1}\) in
    \eqref{eq:fdyn-parallel-round}, we obtain
    \[
    \mathbb E\!\left[\max_{i\in S^{r-1}}K_i^r\right]
    =\mathcal O\!\left(
    \log\!\left(2+s\kappa^2
    \left(\frac{\rho}{\rho-\delta}\right)^{r-1}\right)\right).
    \]
    Since \(\log(\rho/(\rho-\delta))=\mathcal O(\mu/\gamma)\),
    summing over the rounds gives
    \[
    \mathbb E[G_{\rm par}(R)]
    =\mathcal O\!\left(
    R\log(2+s\kappa^2)+\frac{\mu}{\gamma}R^2\right).
    \]
    Substituting the communication bound yields
    \[
    \mathbb E[G_{\rm par}(R)]
    =\mathcal O\!\left(
    \left(\kappa+\frac Ns\right)\log\frac{1}{\epsilon}
    \left[\log(2+s\kappa^2)+\log\frac{1}{\epsilon}\right]\right).
    \]
    For \(0<\epsilon<1\), the bracket is of order
    \(\log(2+s\kappa^2/\epsilon)\), giving the same overall order
    as in part~(a).
\end{itemize}

\end{appendices}


\bibliography{references}

\end{document}